\documentclass{article} 
\usepackage[letterpaper,margin=1.0in]{geometry}
\usepackage{hyperref}
\usepackage{graphicx} 
\usepackage{caption} 
\usepackage{algorithm}
\usepackage{algorithmic}
\usepackage[utf8]{inputenc}
\usepackage{amsmath}
\usepackage{amssymb}
\usepackage{microtype}
\usepackage{multirow, hhline}
\usepackage[table]{xcolor}
\usepackage[para,online,flushleft]{threeparttable}

\usepackage{url}            
\usepackage{booktabs}       
\usepackage{amsfonts}       
\usepackage{bbding}
\usepackage{xcolor}         
\usepackage{graphicx}
\usepackage{enumitem}
\usepackage[numbers,sort]{natbib}
\usepackage{color}
\usepackage{amsmath}\allowdisplaybreaks
\usepackage{amsfonts,bm}

\def\1{\bf{1}}

\newcommand{\Norm}[1]{\left\| #1 \right\|}
\newcommand{\norm}[1]{\left\| #1 \right\|_2}
\def\inner#1#2{\langle #1, #2 \rangle}
\def\Inner#1#2{\left\langle #1, #2 \right\rangle}

\def \mSigma{{\bf \Sigma}}
\def\vzero{{\bf{0}}}

\def\va{{\bf{a}}}

\def\vx{{\bf{x}}}
\def\vy{{\bf{y}}}
\def\vz{{\bf{z}}}

\def\fI{{\mathcal{I}}}

\def\fO{{\mathcal{O}}}

\def\fS{{\mathcal{S}}}

\def\BR{{\mathbb{R}}}

\def\mF {{\bf F}}
\def\mG {{\bf G}}

\def\mJ {{\bf J}}

\usepackage{amsthm}
\usepackage{etoolbox}
\theoremstyle{plain}

\makeatletter
\def\Ddots{\mathinner{\mkern1mu\raise\p@
\vbox{\kern7\p@\hbox{.}}\mkern2mu
\raise4\p@\hbox{.}\mkern2mu\raise7\p@\hbox{.}\mkern1mu}}
\makeatother

\makeatletter
\newcommand*{\rom}[1]{\expandafter\@slowromancap\romannumeral #1@}
\makeatother

\newtheorem{theorem}{Theorem}
\newtheorem{proposition}[theorem]{Proposition}
\newtheorem{lemma}[theorem]{Lemma}
\newtheorem{corollary}[theorem]{Corollary}
\theoremstyle{definition}
\newtheorem{definition}[theorem]{Definition}
\newtheorem{assumption}[theorem]{Assumption}
\theoremstyle{remark}
\newtheorem{remark}[theorem]{Remark}

\def\vx {{{\bf x}}}

\def\vy {{{\bf y}}}
\def\va {{{\bf a}}}
\def\BR{{\mathbb{R}}}

\def\F{{\bf F}}

\def\G{{\bf G}}

\def\I{{\bf I}}
\def \J{{\bf J}}

\def\U{{\bf U}}

\def\V{{\bf V}}

\def\x{{\bf x}}

\def\y{{\bf y}}

\def\z{{\bf z}}
\def\0{{\bf 0}}
\def\1{{\bf 1}}

\def\OM{{\mathcal O}}

\def\RB{{\mathbb R}}

\usepackage{algorithm}
\usepackage{algorithmic}
\usepackage{enumitem}
\usepackage{hhline}
\def\defeq{\overset{\text{def}}{=}}

\def \Floor #1{{\left \lfloor #1  \right \rfloor }}
\usepackage{newfloat}
\title{An Enhanced Levenberg--Marquardt Method via Gram Reduction}

\author{Chengchang Liu\textsuperscript{1}\qquad Luo Luo\textsuperscript{2, 3}\qquad John C.S. Lui\textsuperscript{1}
\\
\vspace{-2mm} \\
\normalsize{
\textsuperscript{1} The Chinese University of Hong Kong
\quad 
\textsuperscript{2} Fudan University}\\
\normalsize{
\textsuperscript{3} Shanghai Key Laboratory for Contemporary Applied Mathematics
}\\
 \vspace{-2mm} \\
\normalsize{ \texttt{7liuchengchang@gmail.com,  luoluo@fudan.edu.cn,
cslui@cse.cuhk.edu.hk
}}
}

\date{}

\begin{document}

\maketitle
\begin{abstract}
  This paper studied the problem of solving the system of nonlinear equations ${\bf F}(\vx)=\vzero$, where $\mF:{\mathbb R}^{d}\to{\mathbb R}^d$. 
  We propose Gram-Reduced Levenberg--Marquardt method which updates the Gram matrix ${\bf J}(\cdot)^\top{\bf J}(\cdot)$ in every~$m$ iterations, where ${\bf J}(\cdot)$ is the Jacobian of ${\bf F}(\cdot)$.
  Our method has  a global convergence guarantee without relying on any step of line-search or solving  sub-problems.   
  We prove our method takes at most $\OM(m^2+m^{-0.5}\epsilon^{-2.5})$ iterations to find an $\epsilon$-stationary point of $\frac{1}{2}\|\F(\cdot)\|^2$, which leads to overall computation cost of~$\OM(d^3\epsilon^{-1}+d^2\epsilon^{-2})$ by taking $m=\Theta(\epsilon^{-1})$. 
  Our results are strictly better than the cost of $\OM(d^3\epsilon^{-2})$ for existing Levenberg--Marquardt methods. 
  We also show the proposed method enjoys local superlinear convergence rate under the non-degenerate assumption.  
  We provide experiments on real-world applications in scientific computing and machine learning to validate the efficiency of the proposed methods.
\end{abstract}

\section{Introduction}
We consider solving the system of  nonlinear equations 
\begin{align}
\label{eq:problem}
    \F(\x)=\0,
\end{align}
where $\x\in\RB^{d}$ and $\F(\x)\triangleq[F_1(\x),\dots,F_d(\x)]^\top\!:\RB^{d}\to\RB^d$ with differentiable $F_i(\vx):\BR^d\to\BR$ for all~$i\in[d]$.
This problem can also be  reformulated by the nonlinear least-square problem:
\begin{align}
\label{eq:problem-2}    \min_{\x\in\RB^d}\phi(\x)\triangleq \frac{1}{2}\|\F(\x)\|^2.
\end{align}
Solving nonlinear equations is one of the most fundamental problems in scientific computing~\cite{nesterov2007modified,yuan2011recent}. 
It has wide applications in machine learning~\cite{defossez2015averaged,botev2017practical,bai2019deep,liu2022loss}, control system~\cite{berthier2021fast}, data assimilation~\cite{tremolet2007model}, and game theory~\cite{frehse1984nonlinear,nourian2013e}.

First-order methods are popular to solve nonlinear equations. Specifically, one can perform the gradient descent on problem (\ref{eq:problem-2}), that is
\begin{align}\label{eq:first}
\begin{split}    
     \x_{t+1} = \x_t-\eta \nabla \phi(\x_t) = \x_t-\eta \J(\x_t)^{\top}\F(\x_t),
\end{split}     
\end{align}
where $\eta>0$ is the step-size, $\phi(\x_t)=\J(\x_t)^{\top}\F(\x_t)$ and~$\J(\cdot)\triangleq[\nabla F_1(\cdot),\dots,\nabla F_d(\cdot)]^{\top}\in\BR^{d\times d}$ is the Jacobian of~$\F(\cdot)$. 
The iteration scheme (\ref{eq:first}) takes $\OM(d^2)$ flops in general. 
By assuming the Lipschitz continuity on $\nabla \phi(\cdot)$, it requires~$\OM(\epsilon^{-2})$ iterations to find an $\epsilon$-stationary point of~$\phi(\cdot)$, then the overall computation cost is $\OM(d^2\epsilon^{-2})$. 
{Gauss--Newton methods \cite{ben1966newton, nocedal1999numerical}} consider the estimation $\nabla^2\phi(\cdot)\approx \J(\cdot)^\top\J(\cdot)$ and establish the iteration
\begin{align*}
     \x_{t+1} = \x_t- \big(\J(\x_t)^{\top}\J(\x_t)\big)^\dag\J(\vx_t)^\top\F(\x_t)
\end{align*}
with local superlinear convergence, while such needs  the flops of $\fO(d^3)$ to compute the matrix (pseudo) inverse.
Recently, \citet{liu2022quasi} introduced Broyden family updates~\cite{rodomanov2021greedy,lin2022explicit} to estimate $\J(\cdot)^{\top}\J(\cdot)$, which leads to the Quasi-Newton methods with local superlinear convergence and reduces the computation cost in per iteration to~$\fO(d^2)$ flops.
However, both Gauss--Newton method and its quasi-Newton variants lack global convergence guarantees. 

In this paper, we are interested in Levenberg--Marquardt (LM) methods~\cite{levenberg1944method,marquardt1963algorithm}, which iterates according to 
\begin{align}
    \label{eq: LM}
    \x_{t+1} = \x_{t} - \big(\J(\x_t)^{\top}\J(\x_t)+\lambda_t\I\big)^{-1}\J(\x_t)^{\top}\F(\x_t),
\end{align}
where $\lambda_t>0$ is the regularization parameter.
This method globalizes the Gauss--Newton iteration, and also maintains local superlinear convergence by properly choosing the regularization term $\lambda_t$~\cite{yamashita2001rate, fan2005quadratic}. 
However, the iteration scheme \eqref{eq: LM} takes~$\OM(d^3)$ flops since it contains matrix inversion like Gauss--Newton update. 
Moreover, the upper bound on iteration numbers for existing global convergent LM methods are no better than $\OM(\epsilon^{-2})$~\cite{ueda2010global,zhao2016global,huang2018global,bergou2020convergence,tran2020stochastic}, which makes the total computation cost be $\OM(d^3\epsilon^{-2})$. 
Besides, these methods require the line-search step or solving trust-region sub-problem~\cite{yuan1994trust} to achieve a descent direction, which leads their implementation be complicated.
Recently, \citet{mishchenko2021regularized} proposed an adaptive LM method without any line-search step nor sub-problem solver, which determine the regularization term by
\begin{align}
\label{eq: regu}
\lambda_t\propto\sqrt{\|\J(\x_t)^{\top}\F(\x_t)\|}\,.
\end{align} 
This method can find an $\epsilon$-stationary within at most $\tilde{\OM}(\epsilon^{-2.5})$ iterations under the cubic-growth condition,
which is worse than line-search or trust-region based LM methods. 
Furthermore, its overall computation cost~$\fO(d^3\epsilon^{-2.5})$ is worse than classical LM methods, and the existence of its local superlinear convergence is unknown. 
Build upon this, there raises a natural question that: 
\textit{
Can we design a computation efficient and globally convergent LM method with local superlinear convergence?
}

We give an affirmative answer to the above question by proposing the Gram-Reduced Levenberg--Marquardt (GRLM) method. 
Instead of existing LM methods that compute the Gram matrix $\J(\x_t)^{\top}\J(\x_t)$ at every iteration, our method only updates $\J(\cdot)^\top\J(\cdot)$ at the snapshot point and reuse it in next $m$ iterations. 
The update of our method can be formulated by
\begin{align}
\label{eq:lazy-lm}
\begin{split}
 \!\!\x_{t+1}\!= \!\x_{t}\! -\!  \Big(\!\J\big(\vx_{\pi(t)}\big)\!^{\top}\J\big(\vx_{\pi(t)}\big)\!+\!\lambda_t\I\Big)^{-1}\!\J(\x_t)\!^{\top}\!\F(\x_t),
\end{split}
\end{align}
where $\pi(t)\triangleq m\Floor{t/m}$ and $\lambda_t$ is chosen according to equation \eqref{eq: regu}.
Note that the scheme \eqref{eq:lazy-lm} only needs to compute the matrix $\J(\vx_{\pi(t)})^{\top}\J(\vx_{\pi(t)})$ in the case of $t\equiv 0~({\rm mod}~m)$, which significantly reduces the time of computing the Gram matrix and leads to a better algorithmic complexity.
For the global behavior, our method takes at most overall computation cost of $\OM(d^3\epsilon^{-1}+d^2\epsilon^{-2})$ to find an $\epsilon$-stationary point under the cubic-growth condition, which improves the results of $\OM(d^3\epsilon^{-2})$ of existing LM methods. 
For the local behavior, our method has the  superlinear rate when the solution has non-degenerate Jacobian, which goes beyond the theory of \citet{mishchenko2021regularized}. 
We summarize the main theoretical results of GRLM and related work in Table~\ref{tbl:compare}.

\paragraph{Paper Organization}
The remainder of the paper is organized as follows. 
In Section~\ref{sec:pre} and \ref{sec:related}, we introduce preliminaries and related work respectively. 
In Section~\ref{sec:lazyLM}, we provide the Algorithm details and convergence analysis of GRLM. 
In Section~\ref{sec:exp}, we validate our methods by numerical experiments. 
We conclude our work and discuss the future directions in Section~\ref{sec:conclu}.

\begin{table*}[!t]
	\centering
	\caption{We summarize the globally convergent algorithms for solving the nonlinear equations by their total computation cost for finding the $\epsilon$-stationary point of $\phi(\cdot)$ ({\bf Globally Computation Cost}), if they enjoy local superlinear rates ({\bf Local Superlinear}), and if they are line-search free ({\bf Line-search Free}).}\label{tbl:compare}
	\begin{threeparttable}
\footnotesize\setlength\tabcolsep{5.pt}
\begin{tabular}{ c  c  c  c c}
			\toprule[.1em]
			 \begin{tabular}{c}\bf Method \end{tabular} &  \begin{tabular}{c}\bf Globally \\ \bf Computation Cost  
        			 \end{tabular} &   \begin{tabular}{c}\bf Local \\ \bf superlinear  \end{tabular}  &  \begin{tabular}{c}\bf Line-search\\\bf Free  \end{tabular}& \bf References \\
			\midrule
         \begin{tabular}{c}
      Gradient Descent 
         \end{tabular} 
    & $\OM(d^2\epsilon^{-2})$ & {\XSolidBrush} &{\Checkmark}&  {\citet{nesterov2018lectures}} \\\cmidrule{1-1}
\begin{tabular}{c}{Levernberg--Marquardt (v1)} \end{tabular} & ${\OM(d^3\epsilon^{-2}}\log(\epsilon^{-1}))$ &{\Checkmark} &{\XSolidBrush} & \citet{bergou2020convergence} \\	\cmidrule{1-1}     
\begin{tabular}{c}{Levernberg--Marquardt (v2)} \end{tabular} & ${\OM(d^3\epsilon^{-2.5}\log(\epsilon^{-1}))}$\tnote{$\dag$} &{\XSolidBrush}\tnote{$\ddag$} &{\Checkmark}& \citet{mishchenko2021regularized} \\	    
\midrule
	 \begin{tabular}{c}{GRLM (ours)} \end{tabular} & ${\OM(d^2\epsilon^{-2}+d^3\epsilon^{-1})}$ &{\Checkmark} &{\Checkmark}& Algorithm~\ref{alg:LLM} \\			
			\bottomrule		
	\end{tabular} 
	\begin{tablenotes}
		{\scriptsize     
  		\item [{$\dag$}] This results can be improved to $\OM(d^3\epsilon^{-2.5})$ by using our analysis framework. We present a discussion in Remark~\ref{rmk:discuss}.\\
    \item [{$\ddag$}] Our analysis can also provide the local superlinear convergence rate for Levernberg--Marquardt (v2), which is not proved in \citet{mishchenko2021regularized}. 
    We present the result in Remark~\ref{rmk:local}.
			}
	\end{tablenotes}  	
	\end{threeparttable}
\end{table*}

\section{Preliminaries}
\label{sec:pre}
We use notation $\|\,\cdot\,\|$ to present the Euclidean norm of a vector and spectral norm of a matrix, respectively. 
We use notation $\I$ to present the identity matrix. 
We denote $\sigma_{\min}(\cdot)$ as the smallest singular value of a matrix. 
For the ease of presentation, we denote the Gram matrix of $\F(\cdot)$ by
\begin{align}
\label{eq: H}
    \G(\x)\triangleq \J(\x)^{\top}\J(\x)\succeq \0,
\end{align}
where $\mJ(\vx)\in\BR^{d\times d}$ is the Jacobian of $\mF:\BR^d\to\BR^d$ at point~$\vx\in\BR^d$.
We also define the approximate stationary point of the function $\phi(\cdot)=\frac{1}{2}\|\F(\cdot)\|^2$ as follows.
\begin{definition}
    We say $\x\in\BR^d$ is an $\epsilon$-stationary point of $\phi(\cdot)$ if it satisfies that $\|\J(\x)^{\top}\F(\x)\|\leq \epsilon.$
\end{definition}
We make the following standard assumptions on the Jacobian of $\F(\cdot)$.
\begin{assumption}
\label{ass:lip}
    We assume the Jacobian of $\F(\cdot)$ is bounded and Lipschitz continuous, i.e., we have $\|\J(\x)\|\leq L_1$ and
    \begin{align}
    \label{eq:lip}
       \|\J(\x)-\J(\y)\|\leq L_2\|\x-\y\|
    \end{align}
    hold for all $\x,\y\in\RB^d$ and some constants $L_1, L_2\geq 0$.
\end{assumption}
\begin{proposition}
    If the Jacobian $\J(\cdot)$ satisfies Assumption \ref{ass:lip}, then it holds that
    \begin{align}
    \label{eq:lipF}
    \begin{split}
        \|\F(\x)-\F(\y)\|\leq L_1\|\x-\y\|~~~\text{and}~~~\Norm{\F(\x)-\F(\y)-\J(\y)(\x-\y)}\leq \frac{L_2}{2}\|\x-\y\|^2
    \end{split}
    \end{align}
    for all $\vx,\vy\in\BR^d$.
\end{proposition}
The following proposition means Assumption~\ref{ass:lip} leads to the Lipschitz continuity of $\G(\cdot)$.
\begin{proposition}
\label{prop}
    If the function $\F(\cdot)$ satisfies Assumption~\ref{ass:lip}, we have
    \begin{align}
    \label{eq:prop-G}
   \!\!\!   \|\G(\x)\|\!\leq\! L_1^2~~\text{and}~~\|\G(\y)\!-\!\G(\x)\|\!\leq \!2L_1L_2 \|\x\!-\!\y\|
    \end{align}    
    for all $\x,\y\in\RB^d$.
\end{proposition}

\section{Related Work}
\label{sec:related}
The idea of reusing second-information dates back to 1960s, where \citet{shamanskii1967modification} presented a variant of Newton method that constructs a new Jacobian every $m$ iterations and analyzed its local behavior. 
Later, such idea has been generalized to different type of Newton methods~\cite{fan2013shamanskii,wang2006further,lampariello2001global,adler2020new} and to training the large language models~\cite{liu2023sophia,elbakary2024fed}. 
Particularly, \citet{fan2013shamanskii} proposed Shamanskii Levenberg--Marquardt method in the form of
\begin{align*}
    \x_{t+1} = \x_t - \left(\J(\x_{\pi(t)})^{\top}\J(\x_{\pi(t)})+\mu_t\|\F(\x_{\pi(t)})\|^{\alpha}\I\right)^{-1}\J(\x_{\pi(t)})^{\top}\F(\x_t),
\end{align*}
where $\mu_t>0$ and $\alpha \in [1,2]$,
with a fast local superlinear rate.
However, these methods do not have global theoretical advantage compared to traditional methods.

Recently, \citet{doikov2023second} proposed the regularized Newton method with lazy Hessians for  general minimization problem $\min_{\vx\in\BR^d} \phi(\vx)$, which iterates with
\begin{align}
\label{eq:updatelazy_Hessian}
  \!\!\!\x_{t+1}\!=\!\x_t\!-\!\big(\nabla^2 \phi(\x_{\pi(t)})\!+\!\sqrt{\!c\|\nabla \phi(\x_t)\|}\,\I\big)\!^{-1}\nabla \phi(\x_t)
\end{align}
for some $c>0$. 
However, this method is not suitable to solve our system of nonlinear equations (or its nonlinear least-square formulation) in the following aspects:
\begin{itemize}[leftmargin=0.5cm]
     \item The update \eqref{eq:updatelazy_Hessian} requires accessing the second-order information of the objective. In the view of nonlinear least-square formulation \eqref{eq:problem-2}, we have
    \begin{align*}
        \nabla^2\phi(\x_{\pi(t)}) = \J(\x_{\pi(t)})^{\top}\J(\x_{\pi(t)}) + \sum_{i=1}^d F_i(\x_{\pi(t)})\nabla^2 F_i(\x_{\pi(t)}).
    \end{align*}
    Compared with the cost of LM methods mainly depends on $\mJ(\cdot)$,  accessing the Hessians 
    $\nabla^2 F_i(\cdot),\dots,\nabla^2 F_d(\cdot)$ in above equation may be much more expensive.
    \item The convergence guarantees of \eqref{eq:updatelazy_Hessian} require the convexity of $\phi(\cdot)$ and Lipschitz continuity of $\nabla^2\phi(\cdot)$, while the function $\phi(\cdot)=\frac{1}{2}\norm{\mF(\cdot)}$ in our problem~\eqref{eq:problem-2} is non-convex in general and the popular setting of nonlinear equations (Assumption \ref{ass:lip}) cannot guarantee the Lipschitz continuity on~$\nabla \phi(\cdot)$ nor $\nabla^2 \phi(\cdot)$. 
    \item It is also notable to mention that \citet{doikov2023second} also studied the cubic regularized Newton methods~\cite{nesterov2006cubic} with lazy Hessians for the general non-convex case and \citet{chen2024second} studied the lazy extra Newton methods for the minimax problems, however, these algorithms also requires accessing the Hessian $\nabla^2\phi(\cdot)$ and its analysis is based on the Lipschitz continuity of Hessian. Thus, their theory is not applicable for solving the system of nonlinear equations.
\end{itemize}

\section{The Gram-Reduced Levenberg--Marquardt Method}
\label{sec:lazyLM}

\begin{algorithm*}[t]
\caption{Gram-Reduced Levenberg--Marquardt (GRLM) Method}\label{alg:LLM}
\begin{algorithmic}[1]
\STATE \textbf{Input:} $\x_0$, $c$, $T$, and $m$ \\[0.15cm]
\STATE \textbf{for} $t=0,1,\dots T-1$ \\[0.15cm]
\STATE \quad  $\pi(t)=m\Floor{t/m}$, \quad  $\z_t = \x_{\pi(t)}$\\[0.12cm]
\STATE \quad    $\x_{t+1}=\x_t-\big(\J(\z_t)^{\top}\J(\z_t)+\sqrt{c\|\J(\x_t)^{\top}\nabla f(\x_t)\|}\,\I\big)^{-1}\J(\x_t)^{\top}\F(\x_t)$\\[0.15cm]
\STATE \textbf{end for}\\[0.15cm]
\end{algorithmic}
\end{algorithm*}

We propose the Gram-Reduced Levenberg--Marquardt (GRLM) method in Algorithm~\ref{alg:LLM}.
Our GRLM method takes the Gram matrix $\G(\z_t)=\J(\z_t)^{\top}\J(\z_t)$ as the approximation to $\nabla^2 \phi(\z_t)$, which avoid the exact $\nabla^2\phi(\cdot)$ in regularized Newton with Lazy Hessians~\cite{doikov2023second}.
Additionally, the setting of $\vz_t$ in GRLM does not require accessing $\J(\cdot)^\top\J(\cdot)$ in every iterations, which makes the algorithm be more efficient than existing LM methods.

We then consider the convergence of our GRLM method.
For the ease of presentation, we denote 
\begin{align*}
\lambda_t\triangleq\sqrt{c\|\J(\x_t)^{\top}\F(\x_t)\|}~~~~~\text{and}~~~~~
r_t\triangleq \|\x_{t+1}-\x_t\|.   
\end{align*}
where the constant $c>0$ follows the input of Algorithm~\ref{alg:LLM}.

In contrast to the analysis of \citet{mishchenko2021regularized} that only lower bounds $r_t$ by $\lambda_{t+1}$, the following lemma provides both lower bound and upper bound of $r_t$ by $\lambda_{t}$.
\begin{lemma}
\label{lm:rt_bound_lambda_t}
    Under Assumption~\ref{ass:lip}, Algorithm~\ref{alg:LLM} holds that
   \begin{align}
    \label{eq:boundrlambda}
                \frac{\lambda_t^2}{c(L_1^2+\lambda_t)}\leq r_t\leq  \frac{\lambda_t}{c}.
    \end{align}
\end{lemma}

Compared with the ordinary LM method, our GRLM use~$\G(\z_t)$ instead of $\G(\x_t)$ at each iteration, which leads to some additional error terms.
In the following lemma, we show such terms can be bounded by the distance between the current iteration point to the snapshot point due to the Lipschitz continuity of $\G(\cdot)$ we have proved in Proposition~\ref{prop}. 
\begin{lemma}
\label{lm:pre}
Under Assumption~\ref{ass:lip}, Algorithm~\ref{alg:LLM} holds that
\begin{align*}
        \left\|\G(\x_t)(\x_{t+1}-\x_t)+\J(\x_t)^{\top}\F(\x_t)\right\|\leq  \lambda_t r_t  +  2L_1L_2r_t\|\z_t-\x_t\|,
\end{align*}
    and
\begin{align}
\begin{split}
        \label{eq:lminner}
        \Inner{\J(\x_{t})^{\top}\F(\x_t)+\G(\x_t)(\x_{t+1}-\x_t)}{\x_{t+1}-\x_t}\leq-\lambda_t r_t^2+2L_1L_2r_t^2\|\x_{t}-\z_t\|.
\end{split}
\end{align}
\end{lemma}

\subsection{Global Convergence Analysis}
\label{sec:global}

We establish the global convergence for Algorithm \ref{alg:LLM} in this subsection.
We adopt the following cubic-growth assumption which has been well-studied by~\citet{mishchenko2021regularized}.
\begin{assumption}
\label{ass:cubic-growth}
    We assume the function $\F:\BR^d\to\BR^d$ satisfies that
    \begin{align}
    \label{eq:cubic-growth}
      \!  \|\F(\y)\|^2\!\leq\! \|\F(\x)+\J(\x)(\y-\x)\|^2 \!+\! M\|\y-\x\|^3
    \end{align}
    for all $\x,\y\in\RB^d$, where $M>0$ is some constant.
\end{assumption}

The following lemma show that, by properly setting the regularization parameter $c>0$, we can guarantee the descent property of $\|\F(\x_t)\|$ at the snapshot point $\vx_t\in\BR^d$ such that~$t\equiv 0~({\rm mod}~m)$.
\begin{lemma}
\label{lm:descent}
    Under Assumption~\ref{ass:lip} and \ref{ass:cubic-growth}, if we run Algorithm~\ref{alg:LLM} with $c=\max\{4L_1L_2m,M\}$, then it holds that
    \begin{align*}
        \|\F(\x_{(k+1)m})\|^2\leq \|\F(\x_{km})\|^2\leq \cdots \|\F(\x_{0})\|^2
    \end{align*}
    and 
    \begin{align*}
        \|\F(\x_{km})\|^2-\|\F(\x_{(k+1)m})\|^2\geq \sum_{t={km}}^{(k+1)m-1}\frac{r_t^2\lambda_t}{6},
    \end{align*}
    for all $k=0,1,\cdots$.
\end{lemma}

Now, we present the global convergence results of our GRLM (Algorithm~\ref{alg:LLM}). 
\begin{theorem}
\label{thm:LMglobal}
    Following the settings of  Lemma~\ref{lm:descent}, Algorithm~\ref{alg:LLM} takes at most
    \begin{align*}
        T= \OM(m^2 + m^{-0.5}\epsilon^{-2.5})
    \end{align*}
    iterations to find an $\epsilon$-stationary point of $\phi(\cdot)$, i.e., we have
    \begin{align*}
        \min_{t=0,\cdots,T-1} \|\J(\x_t)^{\top}\F(\x_t)\|\leq \epsilon.
    \end{align*}
\end{theorem}

\begin{remark}
\label{rmk:discuss}
    Taking $m=1$, Algorithm~\ref{alg:LLM} reduces to the ordinary Levernberg--Mardquardt method proposed 
 by~\citet{mishchenko2021regularized}.
   Compared with the iteration complexity of ${\OM}(\epsilon^{-2.5}\log(1/\epsilon))$ established in the previous work, 
   Theorem~\ref{alg:LLM} provides an improved iteration complexity of $\OM(\epsilon^{-2.5})$, which removes the logarithmic factor. 
   We achieve this by lower bounding $r_t$ by $\lambda_t$ rather than $\lambda_{t+1}$ in the existing work and using a novel proof framework that divides the iterations according to the value of $\lambda_t$.
\end{remark}
Recall that Algorithm~\ref{alg:LLM} needs to compute a new $\G(\z_t)$ in every $m$ iterations. We use $K$ to denote the numbers of the snapshot points, then Theorem \ref{thm:LMglobal} means
\begin{align*}
    K = \left\lceil\frac{T}{m}\right\rceil = \OM(m+m^{-1.5}\epsilon^{-2.5}).
\end{align*}
For reusing the Gram matrix $\mG(\vz_t)$ and reducing the computation cost, we can implement Algorithm~\ref{alg:LLM} by performing  singular value decomposition on the Jacobian  
at the snapshot point, which generally takes $\OM(d^3)$ flops and results
\begin{align*}
    \J(\z_t) = \U(\z_t)\mSigma(\z_t)\V(\z_t)^{\top},
\end{align*}
where $\U(\z_t),\V(\z_t)\in\RB^{d\times d}$ are  orthogonal and $\mSigma(\z_t)$ is diagonal. 
Based on the SVD of $\mJ(\vz_t)$, we can update $\x_t$ according to
\begin{align*}
    \x_{t+1} &= \x_t- \left(\G(\z_t)+\lambda_t\I\right)^{-1}\J(\x_t)^{\top}\F(\x_t)\\
    &=\x_t-\V(\z_t)\big((\mSigma(\z_t))^2+\lambda_t\I\big)^{-1}\V(\z_t)^{\top}\J(\x_t)^{\top}\F(\x_t),
\end{align*}
which can be done in $\OM(d^2)$ flops.

Thus, the total computation cost of Algorithm~\ref{alg:LLM} is
\begin{align}
\label{eq:flop}
\begin{split}
    \#\text{flops} =\OM(d^3K+d^2T) =\!\OM(d^3m\!+\!d^3m^{-1.5}\epsilon^{-2.5}\!+\!d^2m^2 \!+\!d^2m^{-0.5}\epsilon^{-2.5}).
\end{split}
\end{align}
Now we desire to select an appropriate $m$ to minimize the overall flops.
We let $m=\OM(d^{\alpha}\epsilon^{-\beta})$ for some $\alpha,\beta\geq0$ and plug it into equation \eqref{eq:flop}, then the total computation cost is 
\begin{align*}
     &\#\text{flops}\\
     &= \OM\big(d^{3+\alpha}\epsilon^{-\beta} + d^{3-1.5\alpha}\epsilon^{-2.5+1.5\beta} + d^{2+2\alpha}\epsilon^{-2\beta}+ d^{2-0.5\alpha}\epsilon^{-2.5+0.5\beta}\big)\\
    & \geq \OM\Big(d^{3-0.25\alpha}\epsilon^{-1.25+0.25\beta}+d^{2+1.25\alpha}\epsilon^{-1.25-0.75\beta}\Big),
\end{align*}
where the inequality is due to the AM–GM and the equality holds when $\alpha=0$ and $\beta=1$.
Therefore, we should take $m=\Theta(\epsilon^{-1})$, which leads to $\#\text{flops} =\OM( d^{3}\epsilon^{-1}+d^2\epsilon^{-2})$.

We formally present above result in the following corollary.
\begin{corollary}
\label{col:best-m}
Following the setting of Lemma~\ref{lm:descent}, running Algorithm~\ref{alg:LLM} with $m=\Theta(\epsilon^{-1})$ can find an~$\epsilon$-stationary point of $\phi(\cdot)$ within $\OM(d^3\epsilon^{-1} + d^2\epsilon^{-2})$ flops.

\end{corollary}
As a comparison, the globally convergent LM method proposed by \citet{mishchenko2021regularized} \mbox{requires} $\OM(d^3\epsilon^{-2.5}\log(1/\epsilon))$ flops to find an~$\epsilon$-stationary point of $\phi(\cdot)$, which is more expensive than our complexity of $\OM( d^{3}\epsilon^{-1}+d^2\epsilon^{-2})$.
Our result is also sharper than the complexity of $\OM(d^3\epsilon^{-2})$ for other LM methods~\cite{ueda2010global,zhao2016global,huang2018global,bergou2020convergence,tran2020stochastic}.

\subsection{Local Convergence Analysis}
We establish the local superlinear rates for Algorithm~\ref{alg:LLM} in this subsection, with the following assumption for our problem.
\begin{assumption}
\label{ass:unique}
We assume that there exists a solution $\vx^*$ for problem \eqref{eq:problem} which has non-degenerate Jacobian, i.e., we have
\begin{align}
\label{eq:unique}
    \F(\x^*)=\0~~~\text{and}~~~ \sigma_{\min}(\J(\x^*))  = \mu
\end{align}
for some $\mu>0$.
\end{assumption}
The following proposition says Assumption \ref{ass:unique} guarantees that the points in the local neighbour of solution $\x^*$ has non-degenerate $\mJ(\cdot)$ and $\G(\cdot)$.
\begin{proposition}[{\cite[Proposition 2.3]{liu2023block}}]\label{prop:local}
    Under Assumption~\ref{ass:lip} and \ref{ass:unique}, for all $\vx\in\BR^d$ satisfying $\|\x-\x^*\|\leq {\mu^2}/{(6L_1L_2)}$, we have
    \begin{align*}
        \sigma_{\min}(\J(\x))\geq \frac{\mu}{\sqrt{2}}\qquad\text{and}\qquad\G(\x)\succeq \frac{\mu^2}{2}\I.
    \end{align*}
\end{proposition}

We first provide the recursion for the distances between the points in local region to the solution.
\begin{lemma}
\label{lm:superlm}
    Under Assumptions~\ref{ass:lip} and \ref{ass:unique}, we suppose some points $\x_t$ and $\vz_t$ generated by Algorithm~\ref{alg:LLM} are sufficient close to the solution $\x^*$ such that 
    \begin{align}
        \label{eq:local-condi}
        \begin{split}
\x_t,\vz_t\in\fS ~\triangleq\!\Big\{ \!\x\!: \!\vx\in\!\BR^d, \|\x\!-\!\x^*\| \!\leq \min\!\big\{\!\frac{\mu^2}{18L_1L_2}, \frac{\mu^4}{64L_1^2c}\big\}\!\Big\},
\end{split}
    \end{align}
  then it holds that
    \begin{align}
    \label{eq:iterlocal}
    \begin{split}
        \|\x_{t+1}-\x^*\|\leq &\alpha_1\|\x_t-\x^*\|^2+\alpha_2\|\x_t-\x^*\|^{1.5} + 2\alpha_1 \|\z_t-\x^*\|\|\x_t-\x^*\|,
    \end{split}
    \end{align}
    where $\alpha_1 \triangleq {L_1L_2}/{\mu^2}$ and  $\alpha_2\triangleq{2L_1c^{0.5}}/{\mu^2}$. Besides, we  have $\x_{t+1}\in\fS$.
\end{lemma}

Now, we show the explicit local superlinear convergence rate of GRLM.

\begin{theorem}
\label{thm:LMlocal}
Under Assumptions~\ref{ass:lip} and \ref{ass:unique}, we run Algorithm~\ref{alg:LLM} with $\x_0$ such that 
    \begin{align}
    \label{eq:initial}
        \|\x_0-\x^*\|\leq \frac{1}{32(\alpha_1+\alpha_2^2)},
    \end{align}
where $\alpha_1$ and $\alpha_2$ follow the definitions in Lemma \ref{lm:superlm}. Then it holds that 
\begin{align*}
     \|\x_t-\x^*\|\leq \frac{1}{2(\alpha_1+\alpha_2^2)} \left(\frac{1}{2}\right)^{2(1+(1+m/2)^{\pi(t)})(1+(t\%m)/2)},
\end{align*}
where $t\% m= t-\pi(t)$.
\end{theorem}
\begin{remark}
\label{rmk:local}
    When $m=1$, Theorem~\ref{thm:LMlocal} provides the superlinear rate of 
    \begin{align*}
        \|\x_t-\x^*\|\leq \frac{1}{2(\alpha_1+\alpha_2^2)}\cdot\left(\frac{1}{2}\right)^{4(1+t/2)},
    \end{align*}
 for the LM method in \citet[Algorithm~2.2]{mishchenko2021regularized}.
\end{remark}

The local convergence behavior has been widely studied for LM methods, while most of work focus on selecting the regularization parameter in the order of $\lambda_t\propto\|\F(\x_t)\|^{\alpha}$~\cite{yamashita2001rate,fan2005quadratic,bergou2020convergence}. 
This paper firstly investigates the local convergence by using $\lambda_t\propto \sqrt{\|\J(\x_t)^{\top}\F(\x_t)\|}$.

The local superlinear rate in Theorem~\ref{thm:LMlocal} is comparable to the rate achieved by regularized Newton with lazy Hessians \cite{doikov2023second}.
Compared with \citet{doikov2023second} where the objective is supposed to be strongly convex, we only assume the Jacobian at the solution is non-degenerate which allows the objective $\phi(\cdot)=\|\mF(\cdot)\|^2$ be non-convex.
Furthermore, we characterize the superlinear convergence of GRLM by the measure of distance to the solution, while the analysis of \citet{doikov2023second} considers the gradient norm of the objective.

The GRLM method (Algorithm~\ref{alg:LLM}) takes 
an average computation cost of $\OM(d^2)$ per iteration when we choose $m=\Omega(d)$, which matches the cost of per iteration in quasi-Newton methods for solving the nonlinear equations~\cite{lin2021explicit,ye2021greedy, liu2022quasi,liu2023block}. 
However, these quasi-Newton methods lack the global convergence guarantees.

\begin{figure}[t]
\centering
\begin{tabular}{ccc}
\includegraphics[scale=0.22]{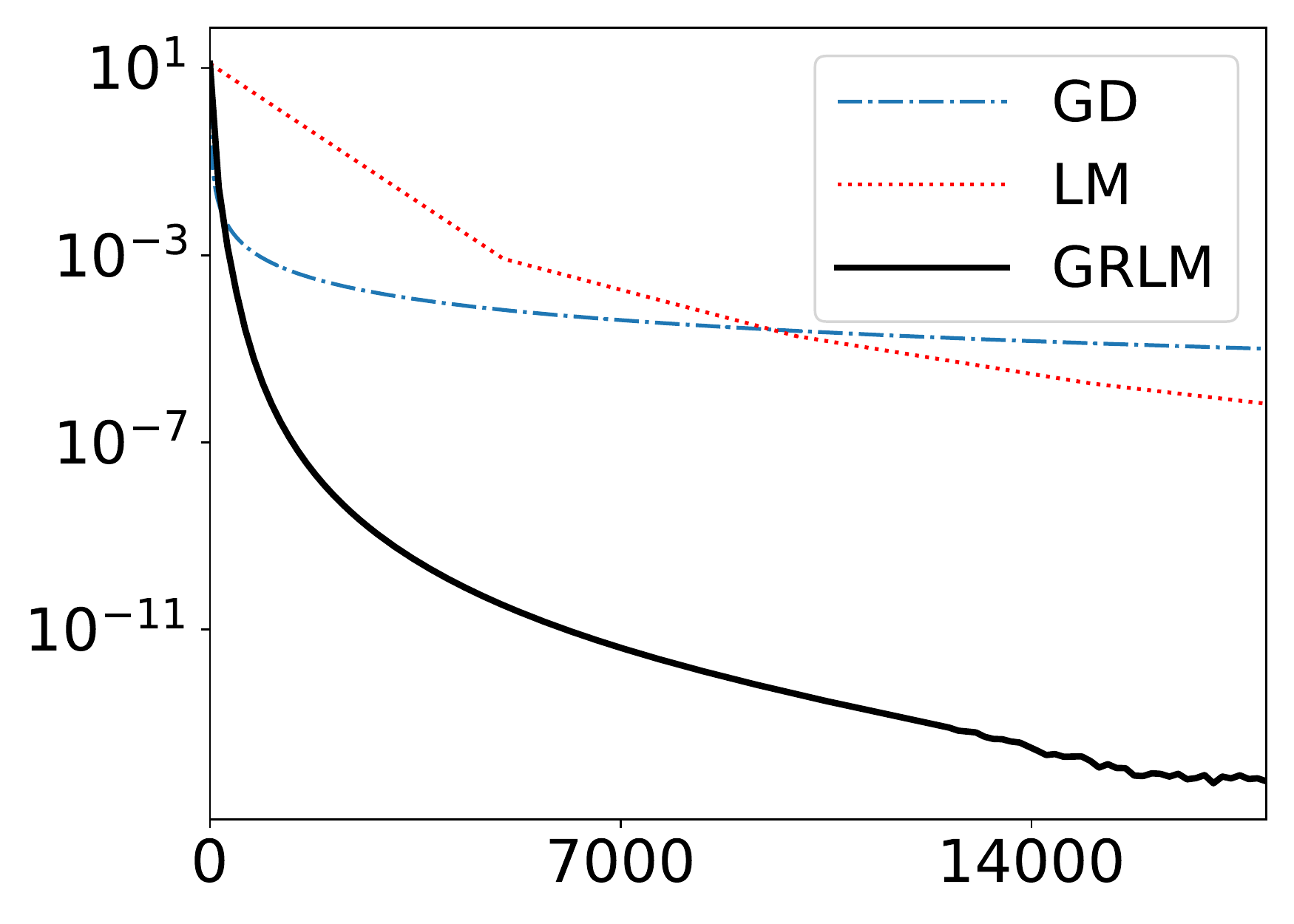} & 
\includegraphics[scale=0.22]{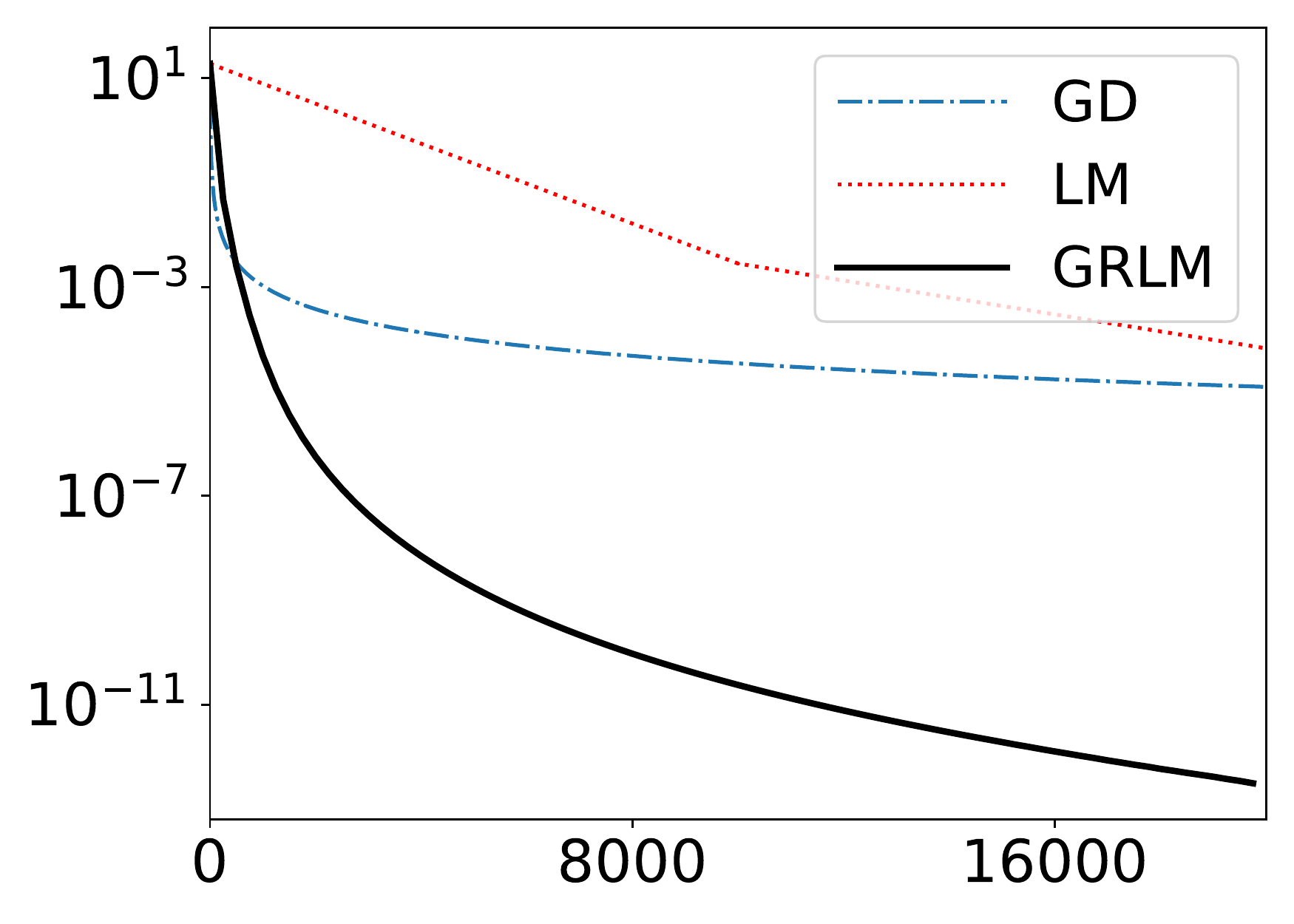} & 
\includegraphics[scale=0.22]{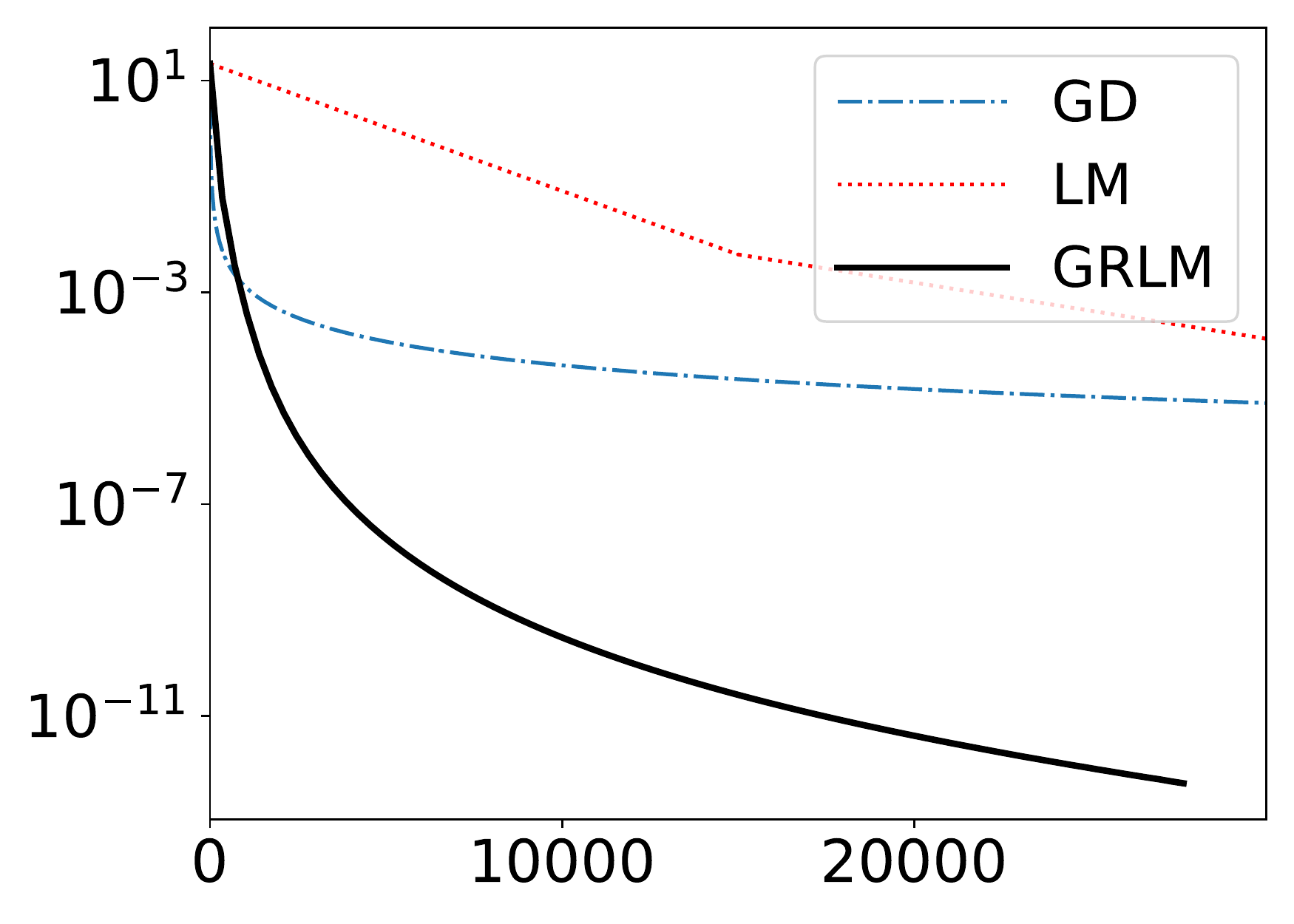} \\
\small (a) $N=100$ ($\#$JV) & \small  (b) $N=200$ ($\#$JV) & \small (c) $N=300$ ($\#$JV)
\\

\includegraphics[scale=0.22]{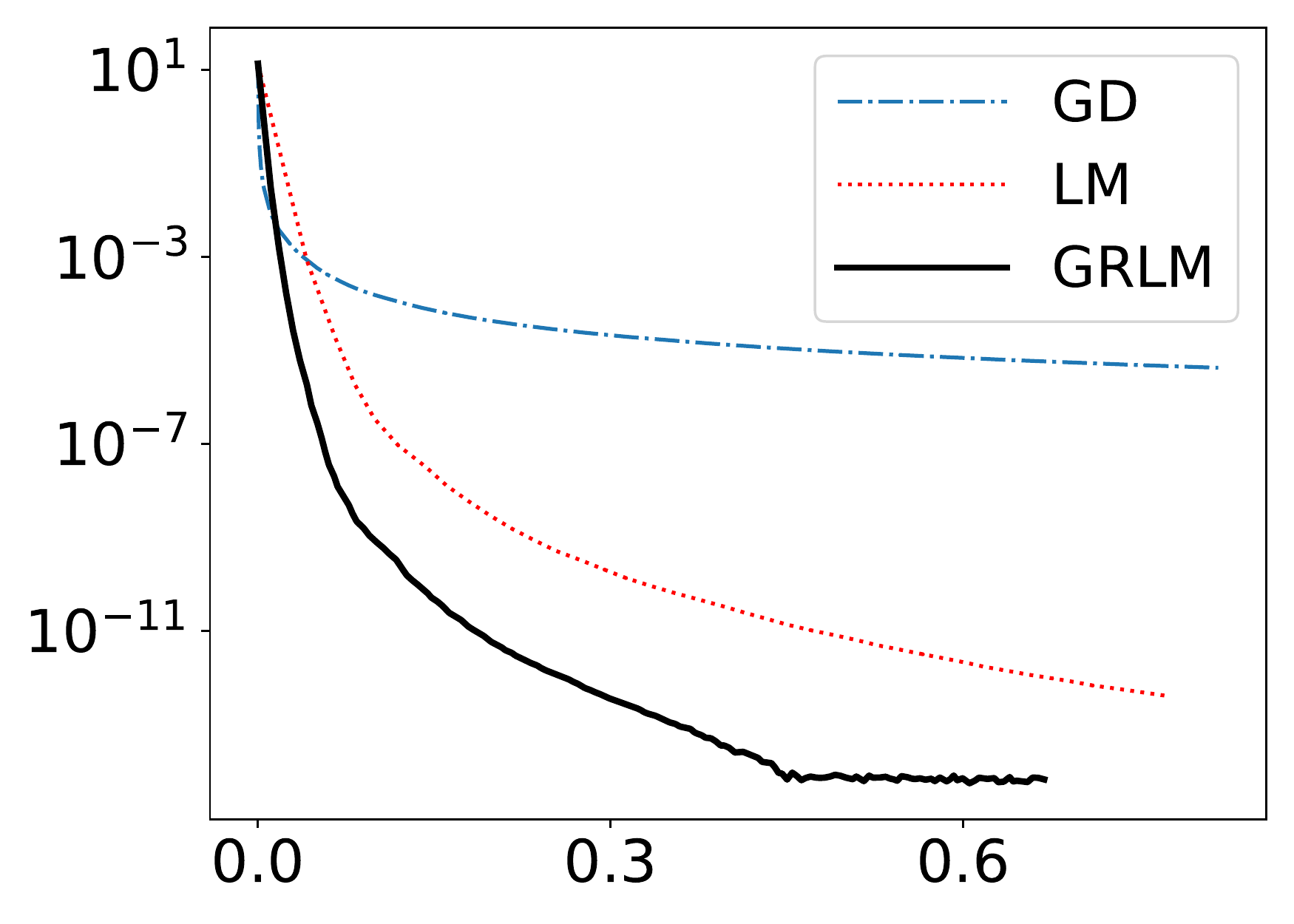} 
& 
\includegraphics[scale=0.22]{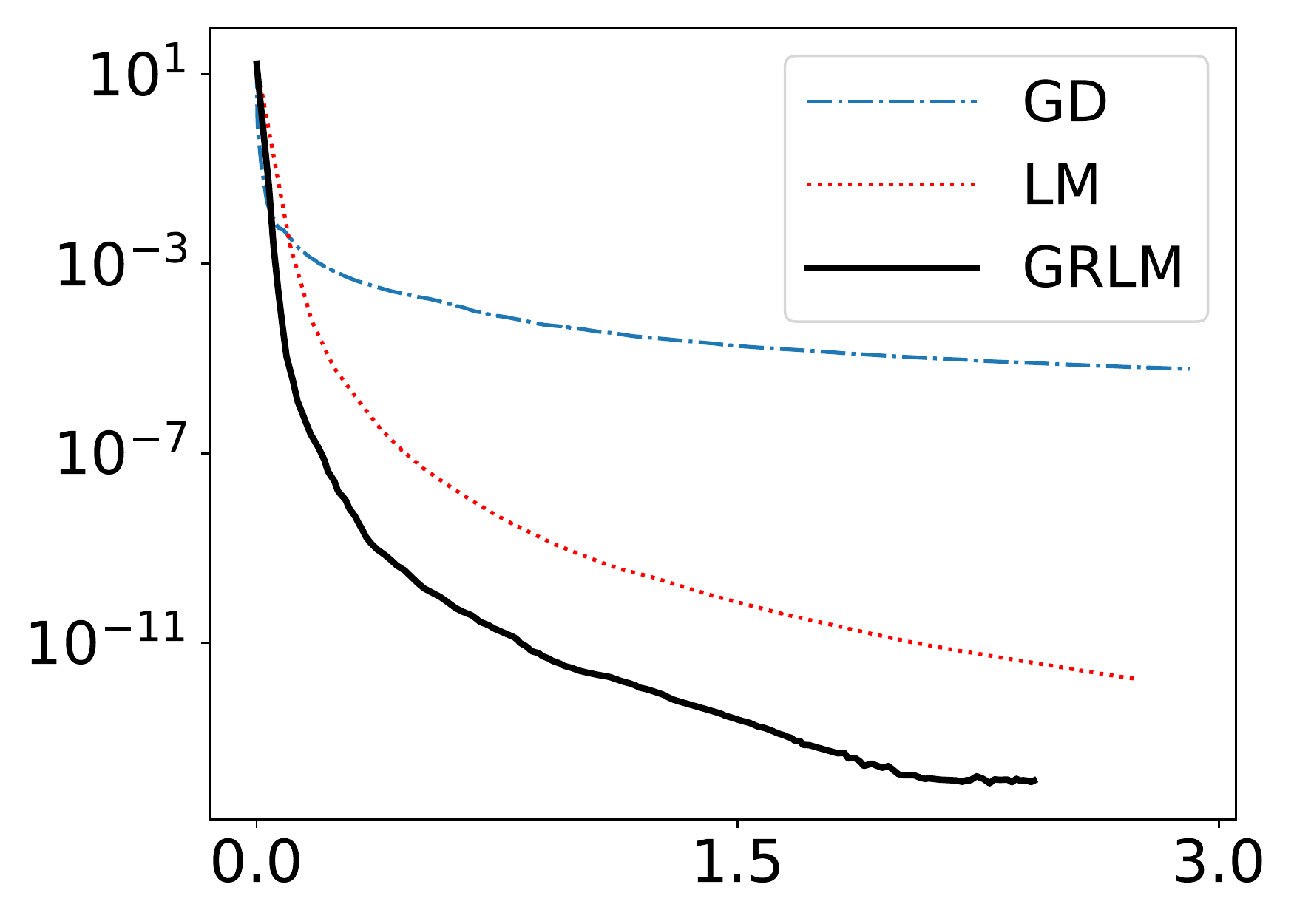}&
\includegraphics[scale=0.22]{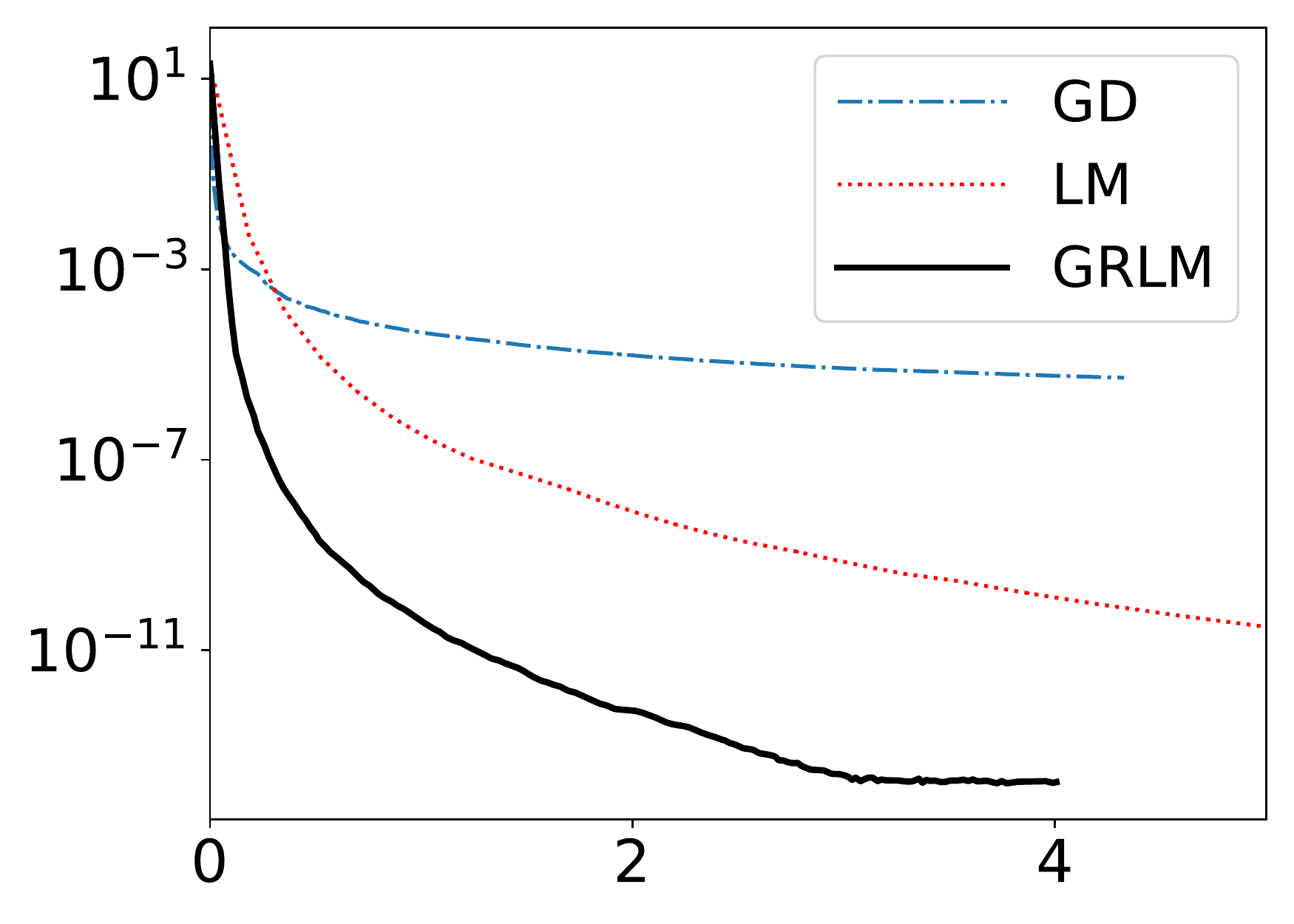} 
\\
 \small  (d) $N=100$ (time) 
 & \small  (e) $N=200$ (time) 
 & \small  (f) $N=300$ (time) 
\\[0.3cm]
\end{tabular}\vskip 0.05cm
\caption{We demonstrate Jacobian-vector products computing times ($\#$JV) and CPU time (second) vs. $\|\J(\x)^{\top}\F(\x)\|_2$ for H-equation with different equation numbers $N$.}\label{fig:compare}\vskip-0.3cm
\end{figure}

\begin{figure}[t]
\centering
\begin{tabular}{ccc}
\includegraphics[scale=0.22]{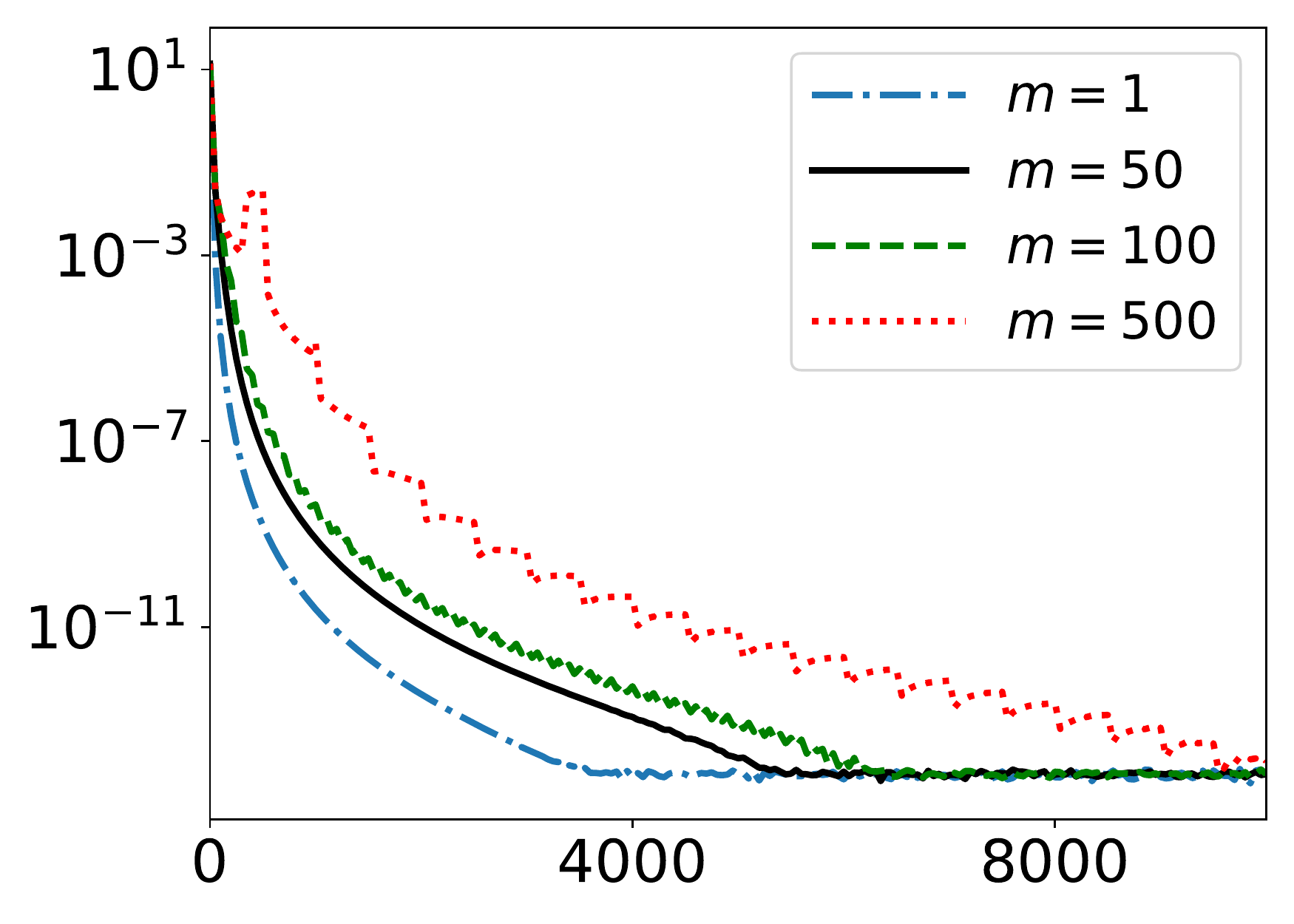} & 
\includegraphics[scale=0.22]{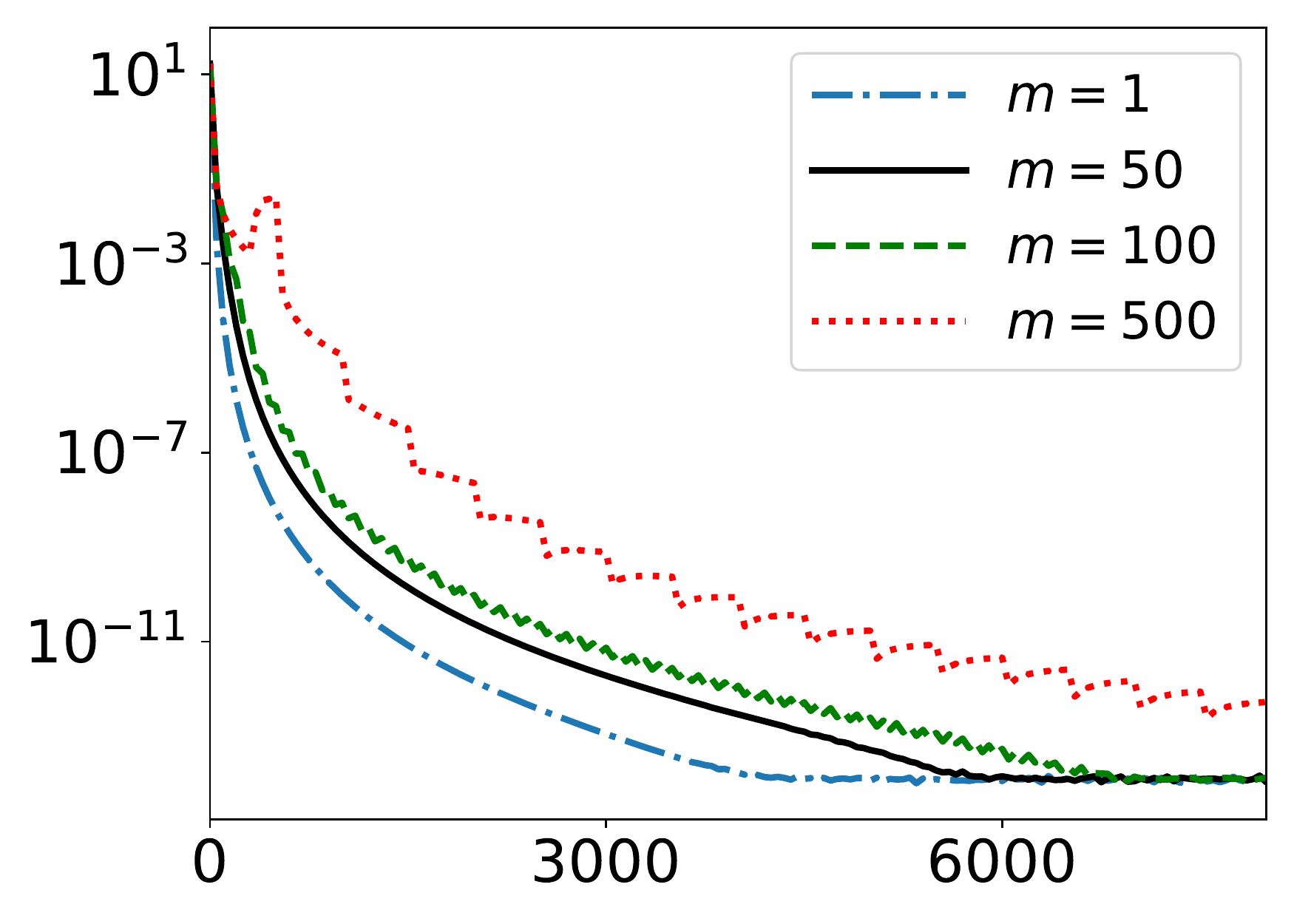} & 
\includegraphics[scale=0.22]{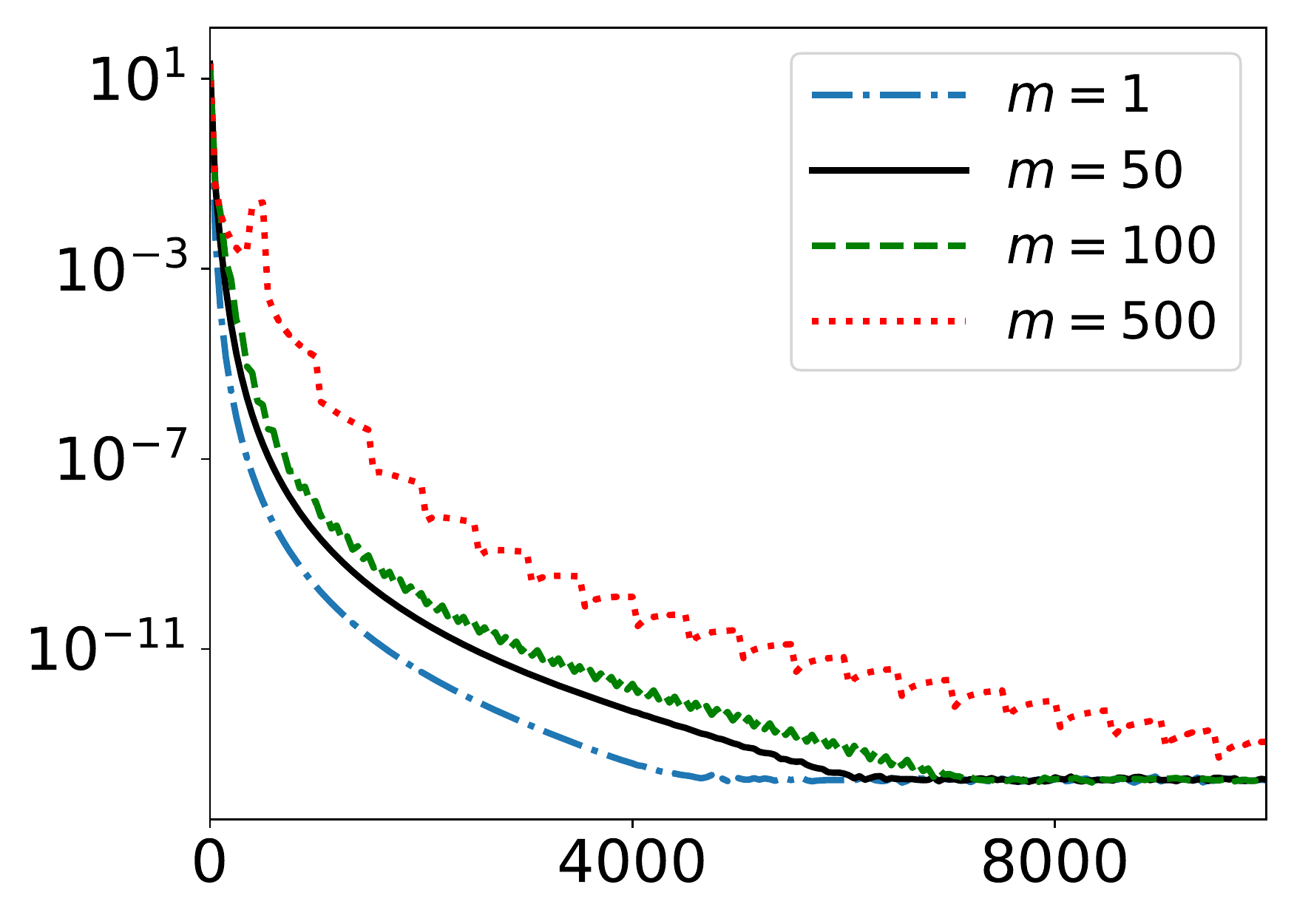} \\
\small (a) $N=100$ (iteration) & \small  (b) $N=200$ (iteration) & \small (c) $N=300$ (iteration)
\\

\includegraphics[scale=0.22]{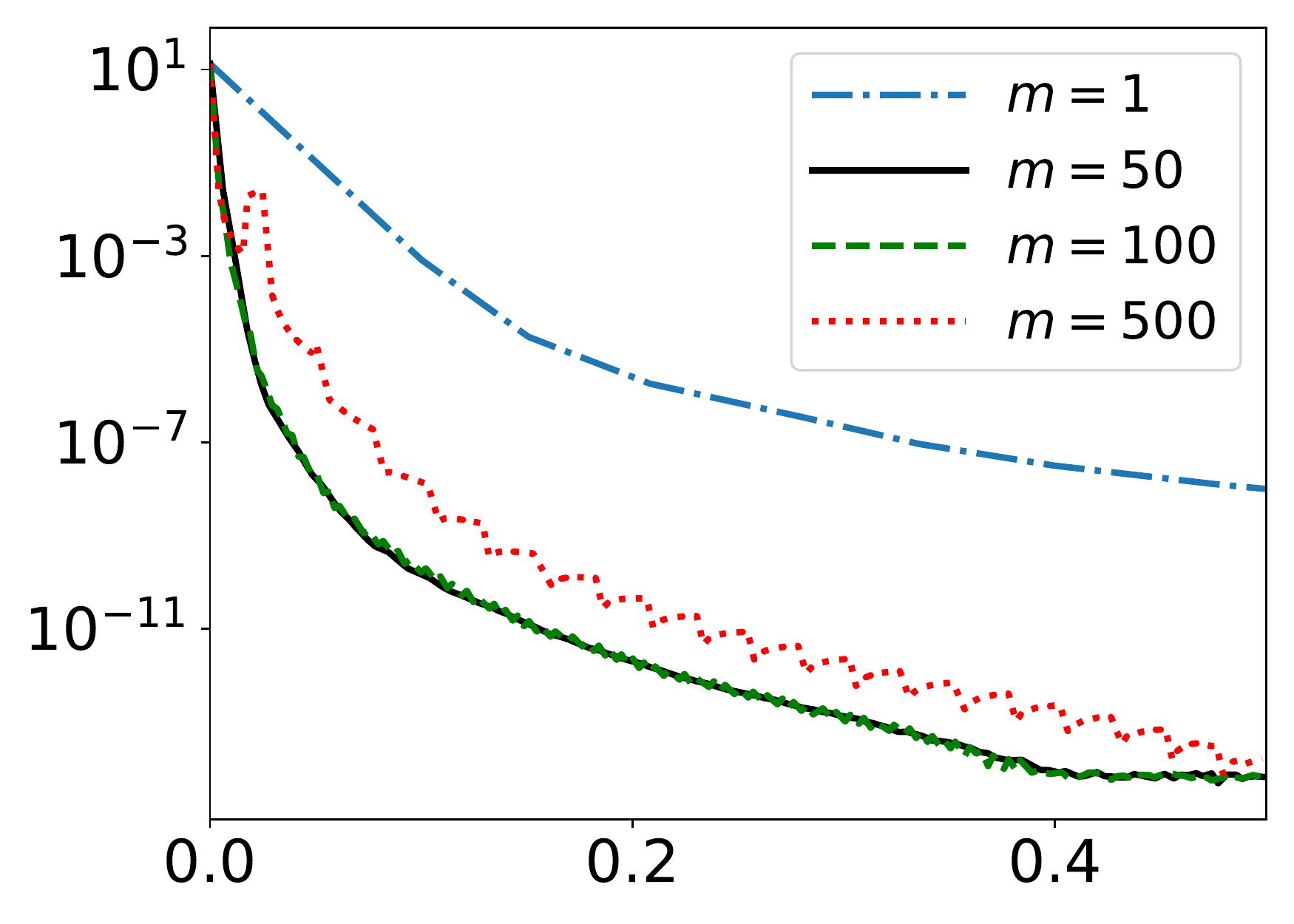} 
& 
\includegraphics[scale=0.22]{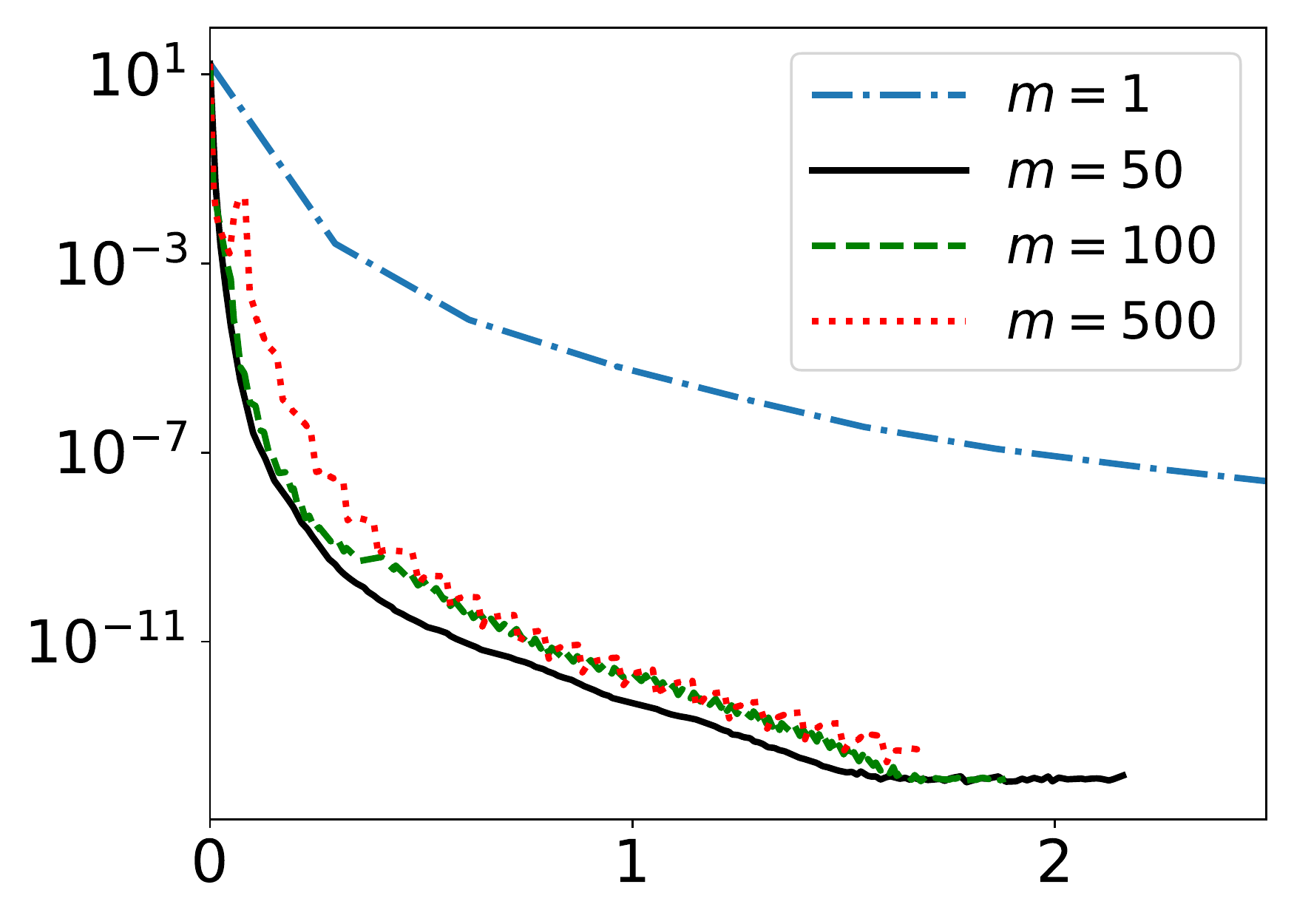}&
\includegraphics[scale=0.22]{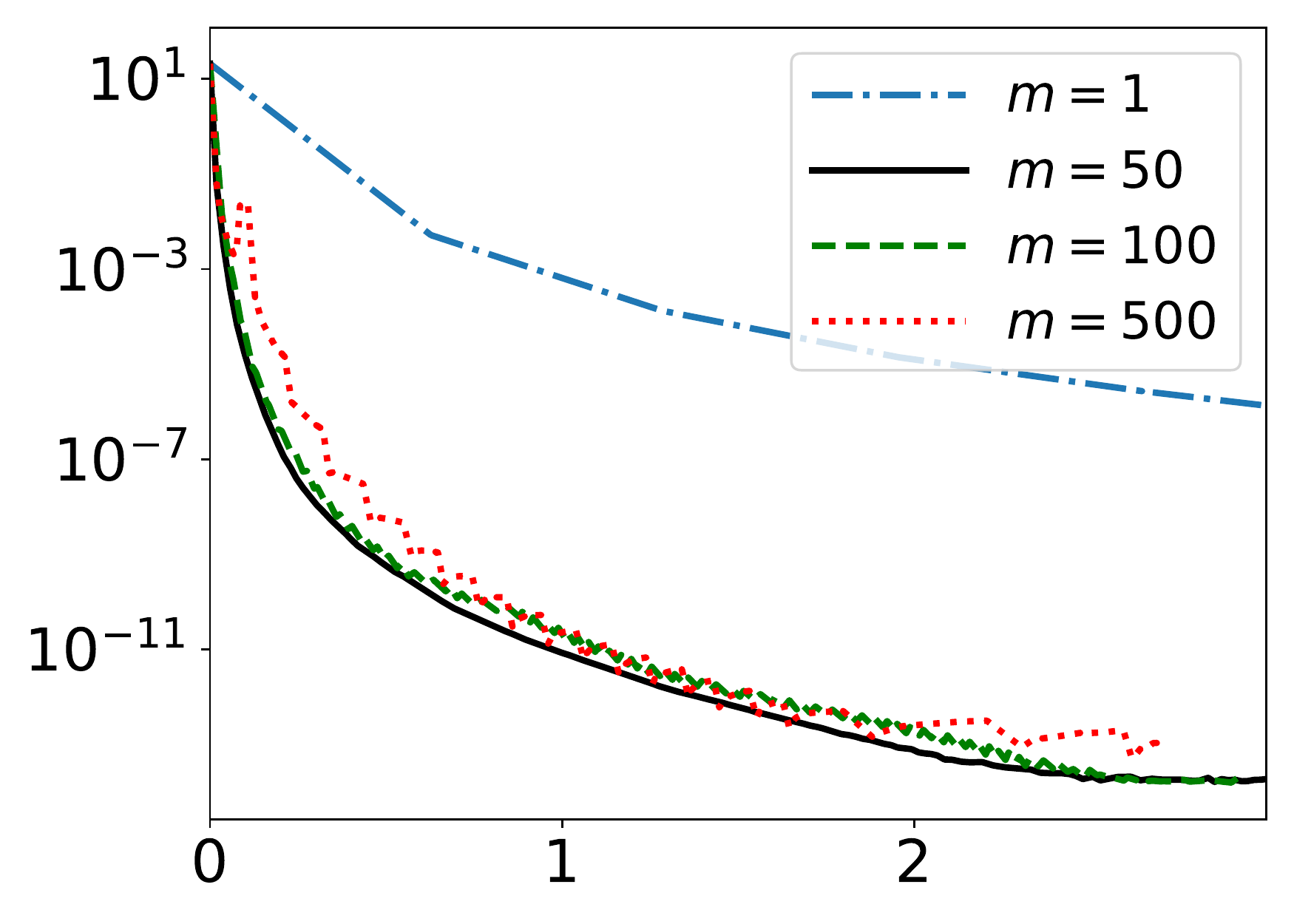} 
\\
 \small  (d) $N=100$ (time) 
 & \small  (e) $N=200$ (time) 
 & \small  (f) $N=300$ (time) 
\\[0.3cm]
\end{tabular}\vskip 0.05cm
\caption{We demonstrate the iteration numbers (iteration) and CPU time (second) vs. $\|\J(\x)^{\top}\F(\x)\|_2$ for H-equation with different equation numbers $N$.}\label{fig:mcompare}\vskip-0.3cm
\end{figure}

\section{Experiment}\label{sec:exp}
In this section, we conduct the experiments on scientific computing and machine learning applications of Chandrasekhar H-equation and non-convex regularized logistic regression to verify our theory.

We choose gradient descent method (GD) and Levernberg--Marquart method with gradient regularization~\cite{mishchenko2021regularized} (LM) as baselines, which are presented in Appendix~\ref{sec:appbase}.
We do not compare our methods with quasi-Newton type methods~\cite{ye2021greedy,lin2021explicit,liu2023block} nor Jacobian-free Newton--Krylov methods~\cite{knoll2004jacobian,ashrafizadeh2015jacobian}, since they do not have global convergence guarantees. 

For all the experiments, we tune the step size $\eta$ for GD from $\{0.1,0.2,\cdots,1\}$. 
We tune the regularized parameter $c$ in LM and GRLM from $\{1,10,100,1000\}$.
Our experiments are conducted on a PC with Apple M1 and all algorithms are implemented in Python 3.8.12.

\subsection{Chandrasekhar H-equation}
\label{sec:H-eq}
Chandrasekhar H-equation plays an important role in scientific computing~\cite{chandrasekhar1960radiative,leggett1976new,kelley1982approximate} and has been well studied in the previous literature~\cite{kelley1995iterative,lin2021explicit,ye2021greedy,liu2023block}. 
It is defined by
\begin{align*}
    F_i(\x)\triangleq x_i -\Big(1-\frac{c}{2N}\sum_{j=1}^{N}\frac{\mu_i x_j}{\mu_i+\mu_j}\Big)^{-1},
\end{align*}
where $\F(\x)=[F_1(\x),\cdots, F_N(\x)]^{\top}\in\RB^{N}$ and  $\x = [x_1,\cdots, x_N]^{\top}\in\RB^{N}$.

We compare GRLM ($m=50$) with the baselines. 
We test the cases $N=100$, $N=200$, and $N=300$. 
In all cases, we set $c=1-10^{-10}$.
We randomize an $\x_0$ as the initial points for all the methods. 
The results of the total number of the Jacobian-vector products ($\#$JV)\footnote{The number of the Jacobian-vector products for computing a full Jacobian is $d$.} computed in the algorithms against $\|\J(\x_t)^{\top}\F(\x_t)\|$ and the running time against $\|\J(\x_t)^{\top}\F(\x_t)\|$ is presented in Figure~\ref{fig:compare}.
We observe that the proposed Gram-reduced Levernberg--Mardquardt method (GRLM) outperforms the baselines in all cases. 

We also provide experiments to study the impact of choosing different $m$ in GRLM (Algorithm~\ref{alg:LLM}).
We choose $m=\{1,50,100,500\}$ for all the cases and present the results in Figure~\ref{fig:mcompare}.

We find that a smaller $m$ leads to a faster iteration in (a), (b), and (c) of Figure~\ref{fig:mcompare}, which matches our theoretical results obtained in Theorem~\ref{thm:LMglobal}.
However, there exists a trade-off between the iteration complexity and the computation cost per-iteration.
For all cases, we find that $m=50$ leads to the smallest CPU time according to the results in (d), (e), (f) of Figure~\ref{fig:mcompare}. 
Such trade-off is consistent to our findings in Corollary~\ref{col:best-m}.

\subsection{Non-Convex Regularized Logistic Regression}
\label{sec:nc-logit}
We further validate the GRLM methods on the non-convex regularized logistic regression model~\cite{antoniadis2011penalized}:
\begin{align*}
    \min_{\x\in\RB^d} f(\x) \triangleq \frac{1}{n}\sum_{i=1}^{n}\ln(1+\exp(-b_i\va_i^{\top}\vx))\! +\! \lambda\sum_{p=1}^d\frac{ x_{(p)}^2}{1+x_{(p)}^2},
\end{align*}
where $x_{(p)}$ is the $p$-th coordinate of $\x\in\RB^d$, $\lambda>0$ is the regularized parameter, $\va_i\in\BR^{d}$ and $b_i\in\{-1,+1\}$ are the feature and the corresponding label of the $i$-th sample. 
This model is corresponding to solve the following nonlinear equations
\begin{align*}
    \F(\x) \triangleq \nabla f(\x).
\end{align*}

\begin{figure}[t]
\centering
\begin{tabular}{ccc}
\includegraphics[scale=0.22]{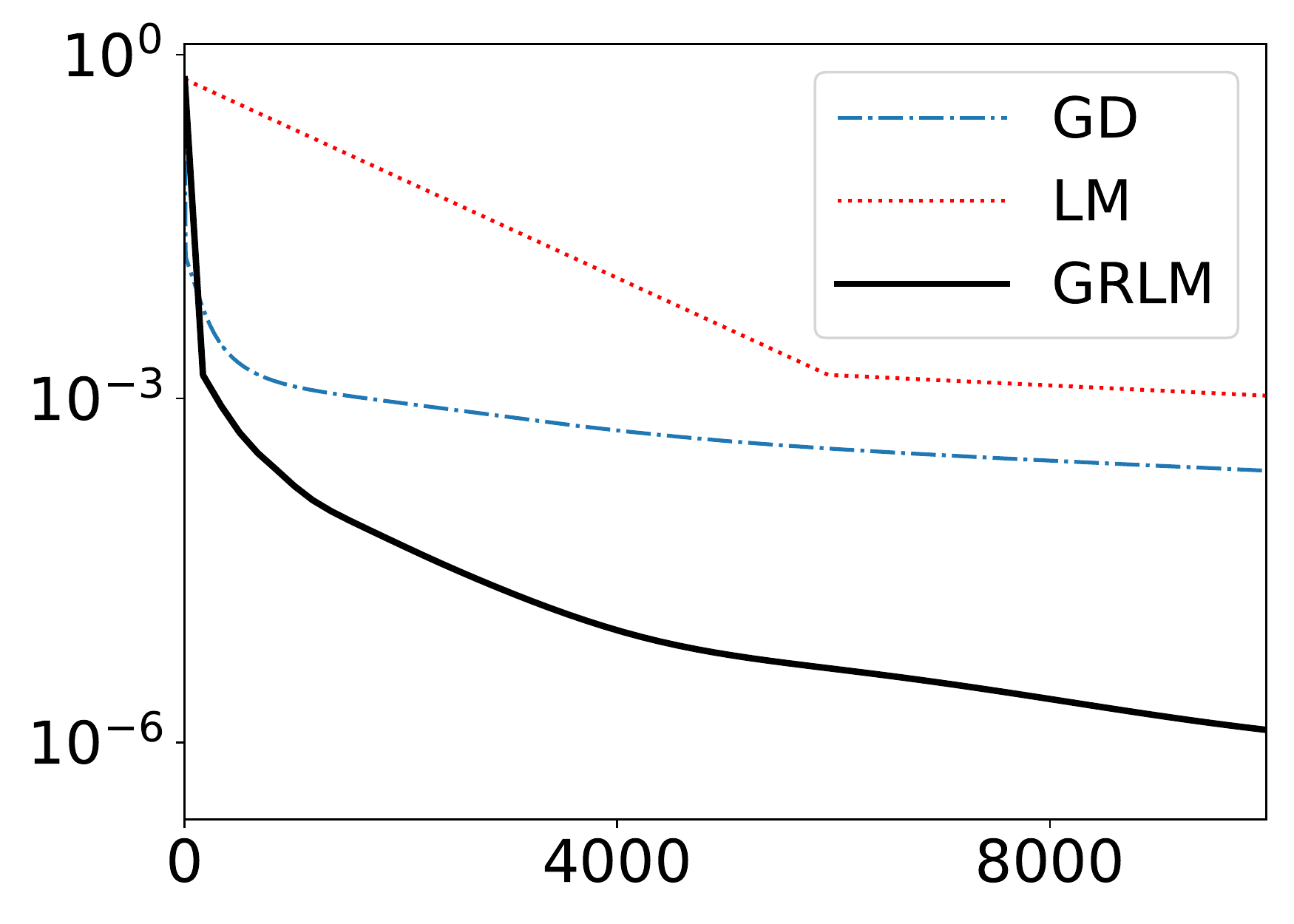} & 
\includegraphics[scale=0.22]{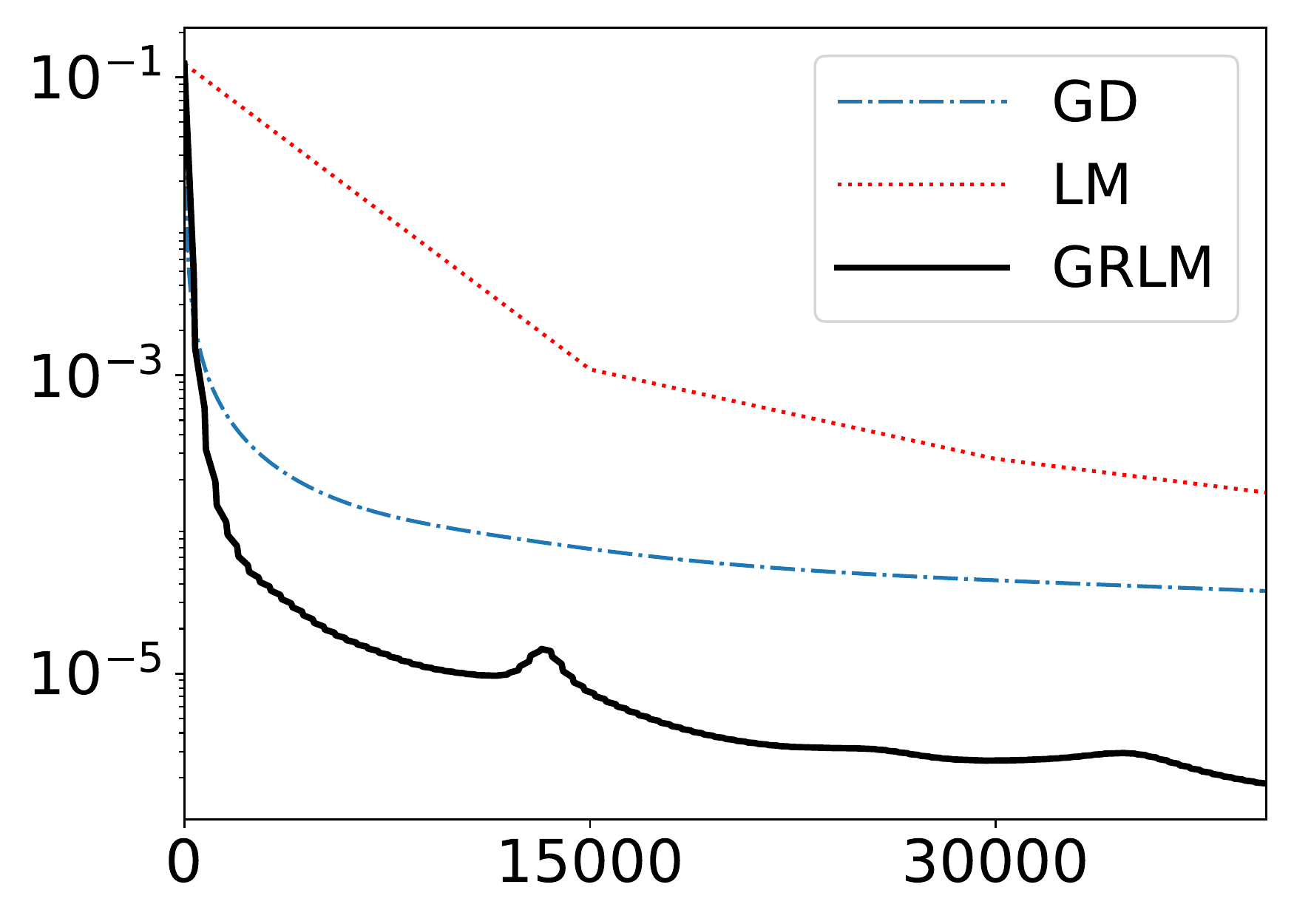} & 
\includegraphics[scale=0.22]{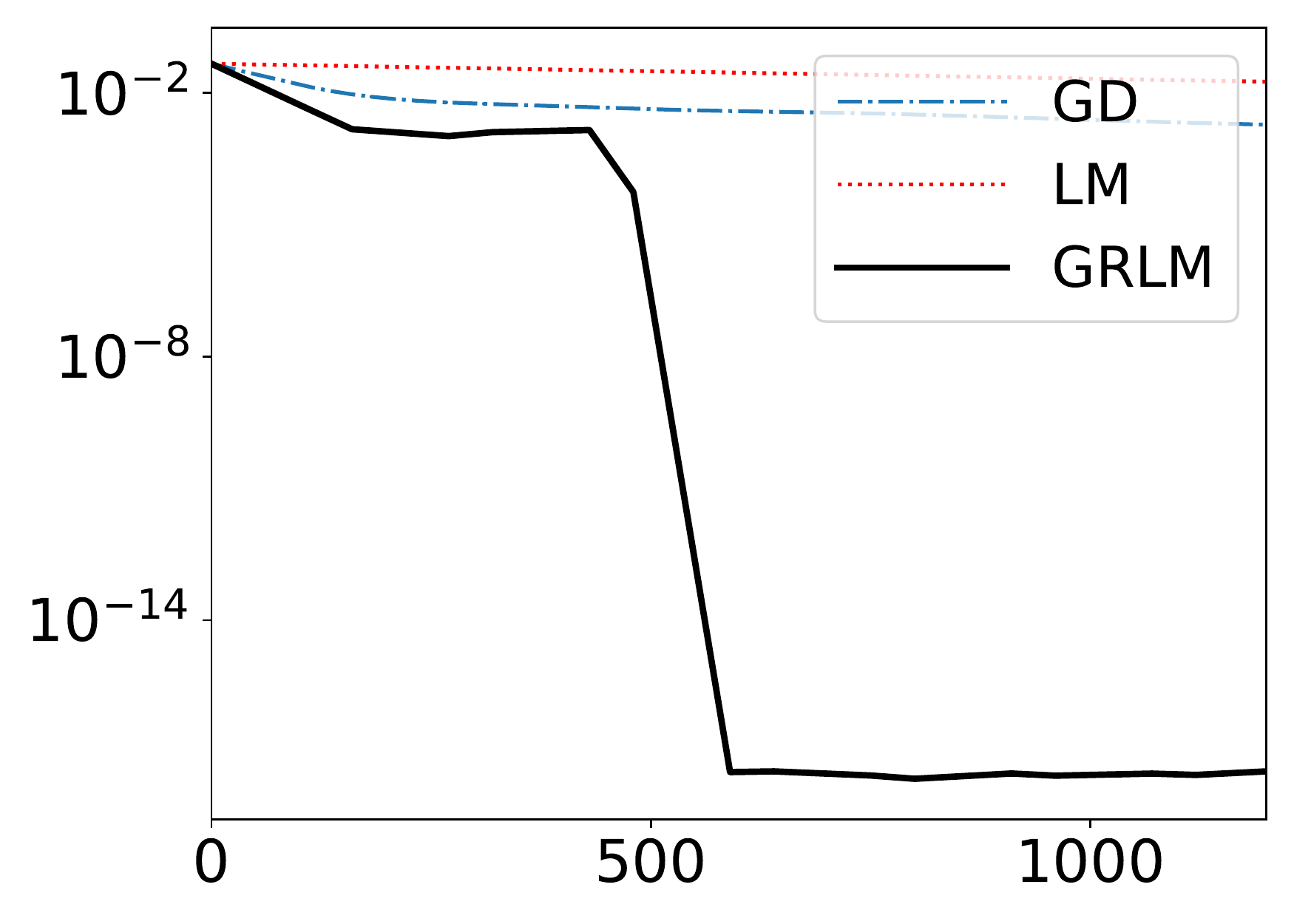} \\
\small (a) ``a1a'' ($\#$JV) & \small  (b) ``w1a'' ($\#$JV) & \small (c) ``splice'' ($\#$JV)
\\

\includegraphics[scale=0.22]{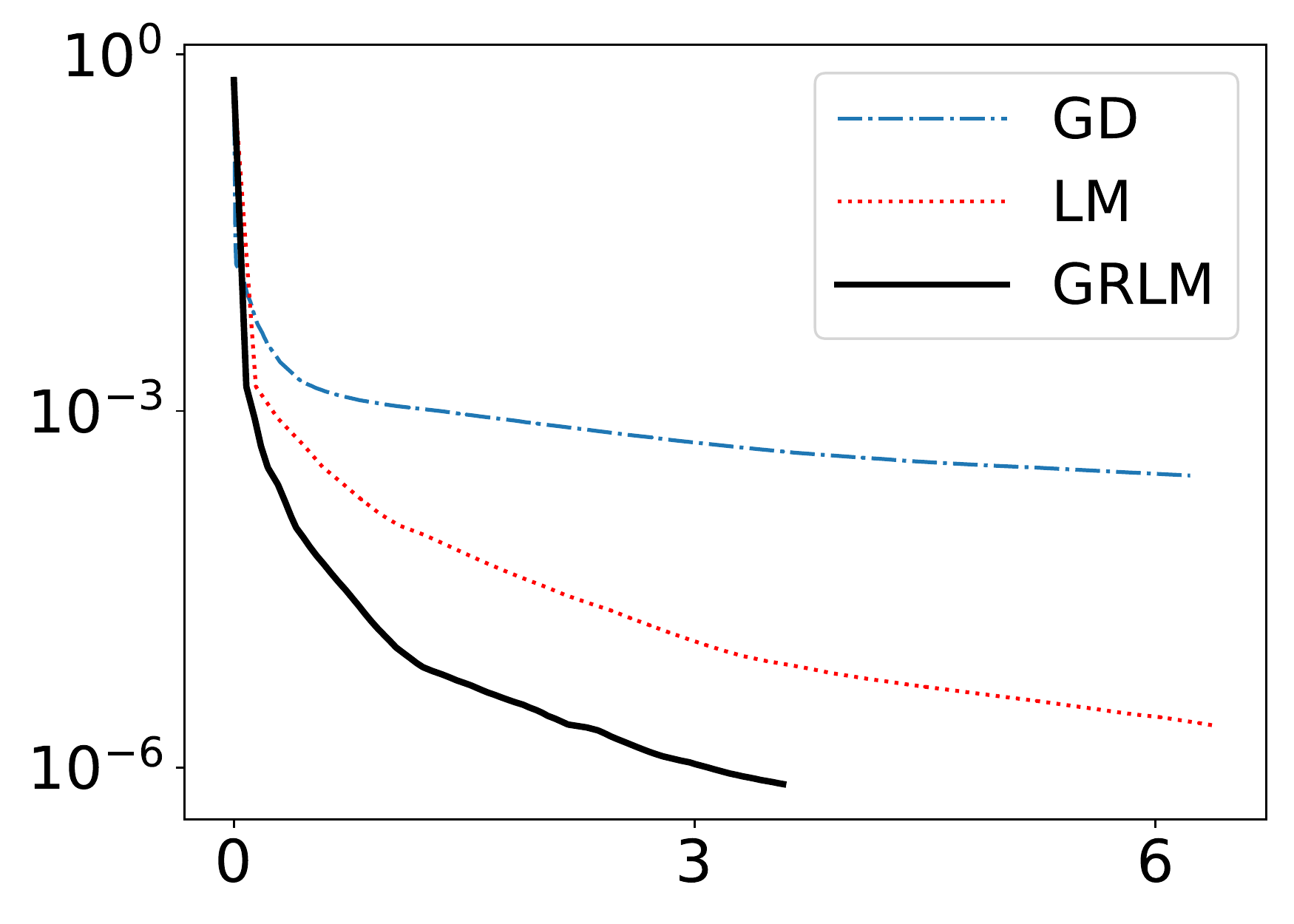} 
& 
\includegraphics[scale=0.22]{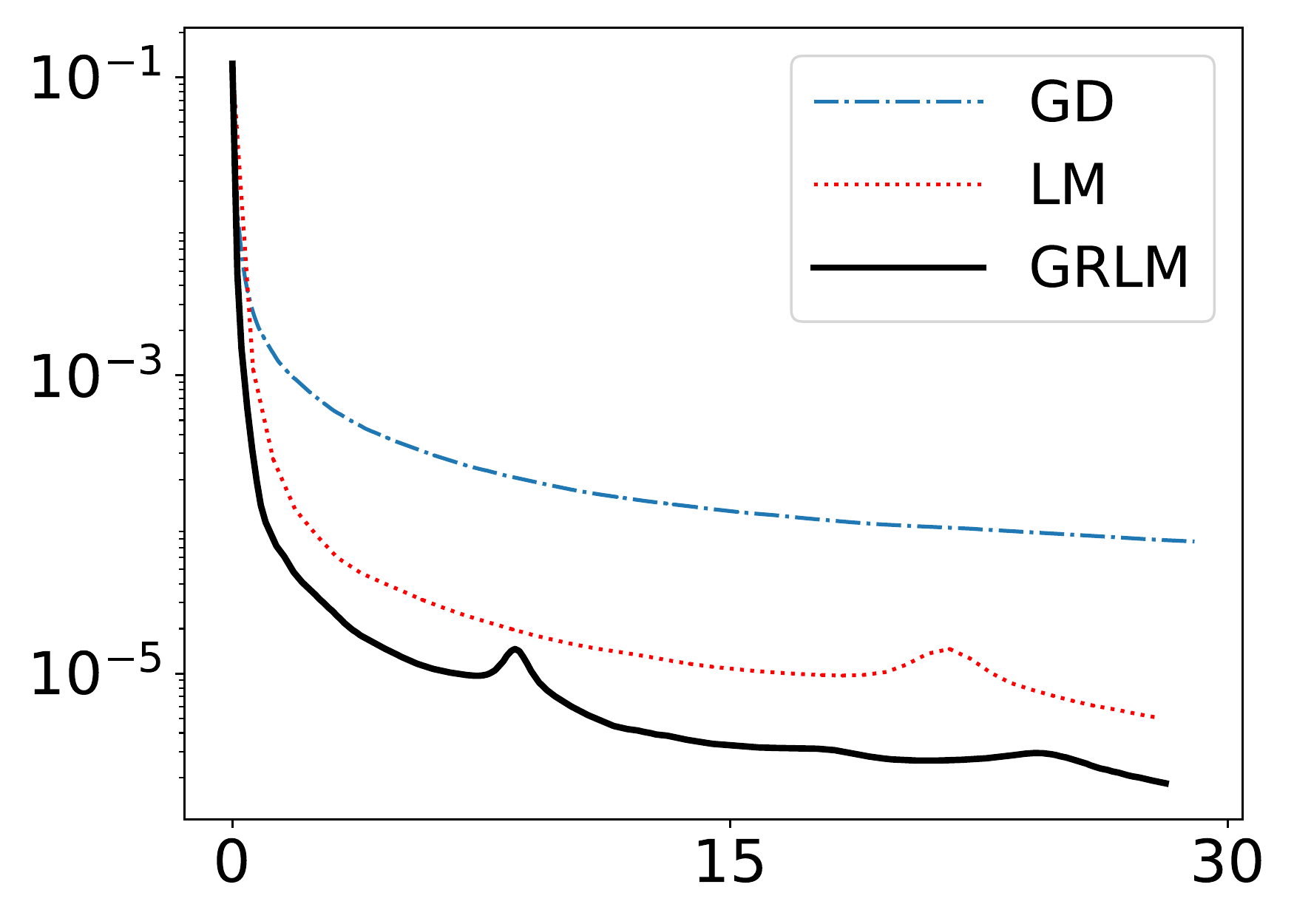}&
\includegraphics[scale=0.22]{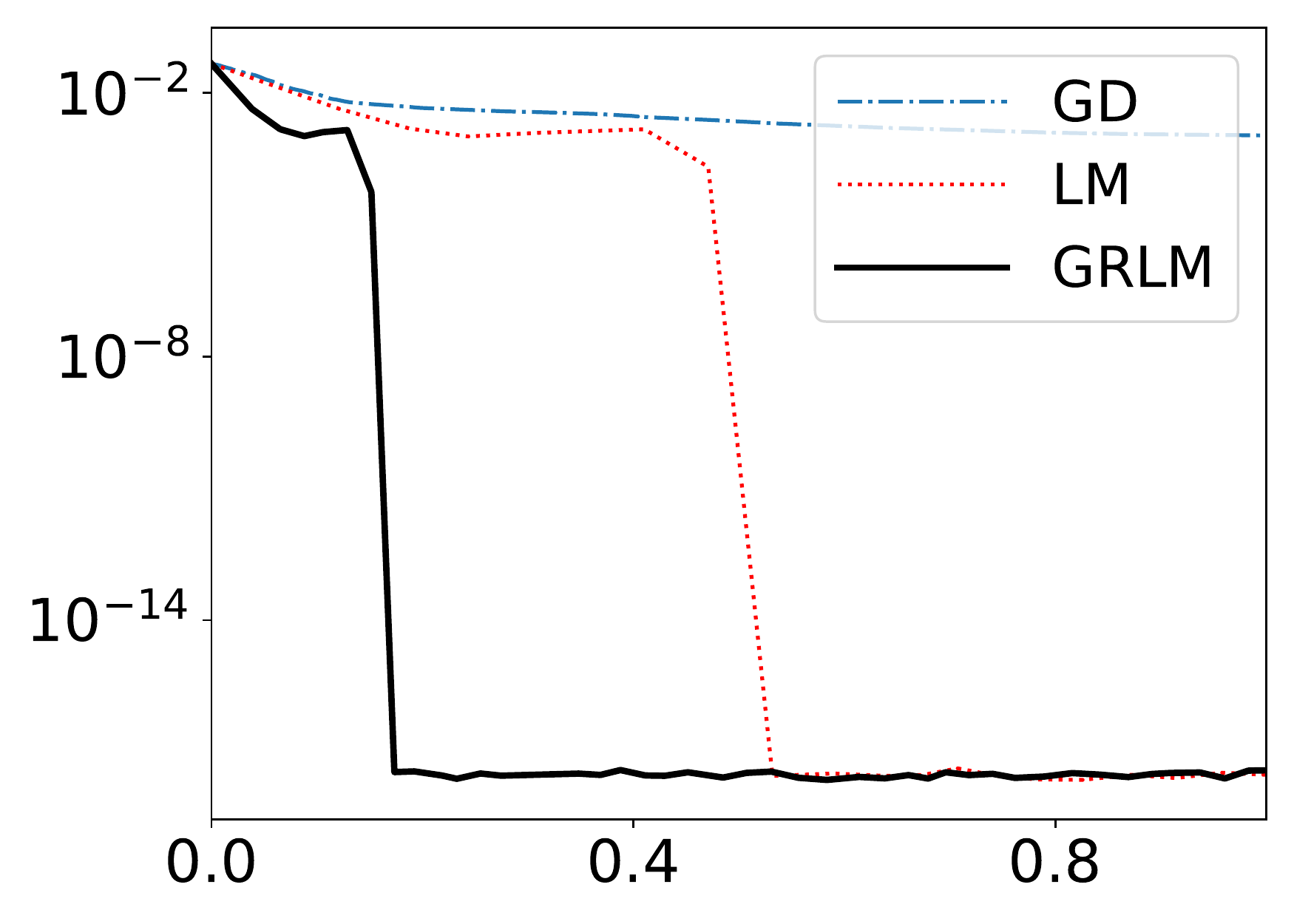} 
\\
 \small  (d) ``a1a'' (time) 
 & \small  (e) ``w1a'' (time) 
 & \small  (f) ``splice'' (time) 
\\[0.3cm]
\end{tabular}\vskip 0.05cm
\caption{We demonstrate Jacobian-vector products computing times ($\#$JV) and CPU time (second) vs. $\|\J(\x)^{\top}\F(\x)\|_2$ for non-convex logistic regression on datasets ``a1a'', ``w1a'', and ``splice''.}\label{fig:compare-NC}\vskip-0.3cm
\end{figure}

We compare GRLM ($m=100$) with the baselines on three real-world datasets: ``a1a'' ($n=1,605$, $d=123$), ``w1a'' ($n=2,477$, $d=300$), and ``splice'' ($n=1,000$, $d=60$). All these datasets can be downloaded from LIBSVM  repository~\cite{CC01a} and are popular for testing the behaviors on the regression model.
We present the results of the number of Jacobian-vector products  ($\#$JV) against $\|\J(\x_t)^{\top}\F(\x_t)\|$ and the running time against $\|\J(\x_t)^{\top}\F(\x_t)\|$ in Figure~\ref{fig:compare-NC}.
GRLM still significantly outperforms the baselines for all the datasets in terms of  both the total numbers of the Jacobian-vector products and the CPU time.

\section{Conclusion}
\label{sec:conclu}
In this paper, we have proposed Gram-Reduced Levenberg--Marquardt method (GRLM) for solving nonlinear equations. 
GRLM enhances the behaviour of existing LM methods by reducing the computation times of the Gram matrices.
It is globally convergent with simple iteration scheme and cheap computation cost.
Furthermore, it exhibits the local superlinear convergence, which cannot be achieved by any of the first-order methods.
To the best of our knowledge, this is the first method for solving nonlinear equations that can achieve the best of the both worlds.

For the future work, it is interesting to design sketched, stochastic, incremental, and distributed variants of GRLM~\cite{yuan2009subspace,yuan2022sketched,chayti2023unified,zhou2024incremental,liu2023communication}. 
It is also possible to incorporate the idea of super universal Newton methods~\cite{doikov2024super} to make GRLM parameter-free.
\bibliographystyle{plainnat}
\bibliography{aaai25.bib}

\appendix

\section{Useful Results for Positive Sequences}
\label{sec:usefulsequece}
We provide useful results for positive sequences which are used in our analysis for self-complement.
\begin{lemma}{\cite[Lemma B.1]{doikov2023second}}
\label{lm:addi-2}
   For any positive sequences $\{r_t\}_{t\geq 1}$, it holds for any~$m\geq 1$:
   \begin{align*}
       \sum_{t=1}^{m-1}\left(\sum_{i=1}^t r_i\right)^3\leq \frac{m^3}{3}\sum_{t=1}^{m-1}r_t^3.
   \end{align*}
\end{lemma}
\begin{lemma}{\cite[Theorem D.3]{doikov2023second}}
\label{lm:addi-1}
For the positive sequences $\{s_t\}_{t=0}^{T}$, if $s_0\leq 1/2^4$ and $s_{t+1}\leq \frac{1}{2}(s_t^2+s_t^{1.5})+s_{\pi(t)}s_t$ holds for all $t\geq 1$, then it holds that
\begin{align*}
    s_{t}\leq \left(\frac{1}{2}\right)^{2(1+(1+m/2)^{\pi(t)})(1+(t\%m)/2)}.
\end{align*}
\end{lemma}

\section{The Proof of Proposition~\ref{prop}}
\label{sec:proof_prop}
\begin{proof}
The norm of the Gram matrix can be bounded by
\begin{align*}
    \|\G(\x)\|\leq \|\J(\x)^{\top}\|\|\J(\x)\|\leq L_1^2.
\end{align*}
The triangle inequality and Assumption~\ref{ass:lip} implies
\begin{align*}
    \|\G(\y)-\G(\x)\|&\leq \|\J(\x)^{\top}(\J(\x)-\J(\y))\|+\|(\J(\x)^{\top}-\J(\y)^{\top})\J(\y)\|\leq 2L_1L_2\|\x-\y\|.
\end{align*}
\end{proof}
\section{The Proof of Lemma~\ref{lm:rt_bound_lambda_t}}
\label{sec:proof_rt_bound}
\begin{proof}
Recall that the update rule of Algorithm~\ref{alg:LLM} implies
\begin{align*}
    \x_{t+1}-\x_t=-\left(\G(\z_t
    )+\lambda_t\I\right)^{-1}\J(\x_t)^{\top}\F(\x_t),
\end{align*}
which means we can upper bound $r_t$ by $\lambda_t$ as follows
\begin{align*}
    r_t &=\|(\G(\z_t
    )+\lambda_t\I)^{-1}\J(\x_t)^{\top}\F(\x_t)\|\\
    &\overset{\eqref{eq: H}}{\leq} \frac{1}{c}\|\lambda_t^{-1}c\J(\x_t)^{\top}\F(\x_t)\| \\
    &= \frac{\lambda_t}{c}.
\end{align*}
We then reorganize the update formula as follows
\begin{align}
\label{eq:reform}
    \G(\z_t)(\x_{t+1}-\x_t) + \lambda_t(\x_{t+1}-\x_t)= -\J(\x_t)^{\top}\F(\x_t) ,
\end{align}
then Assumption \ref{ass:lip} means
\begin{align*}
   \frac{\lambda_t^2}{c}&=\|\J(\x_t)^{\top}\F(\x_t)\|\\
   &\overset{\eqref{eq:reform}}{=} \|\G(\z_t)(\x_t-\x_{t+1}) + \lambda_t (\x_{t}-\x_{t+1})\|\\
   &\leq \|\G(\z_t)(\x_t-\x_{t+1}) \| + \lambda_tr_t \\
   &\leq \|\G(\z_t)\|r_t +\lambda_t r_t\\
   &\overset{\eqref{eq:prop-G}}{ \leq} L_1^2 r_t +\lambda_t r_t,
\end{align*}
which means
\begin{align*}
r_t\geq \frac{\lambda_t^2}{c(L_1^2+\lambda_t)}.
\end{align*}
Thus, we finish the proof of Lemma~\ref{lm:rt_bound_lambda_t}.
\end{proof}

\section{The Proof of Lemma~\ref{lm:pre}}
\label{sec:proof_lemma_pre}
\begin{proof}
We reorganize the update formula of Algorithm~\ref{alg:LLM} as
\begin{align}
\label{eq:re_org}
     \G(\z_t)(\x_{t+1}-\x_t) + \lambda_t(\x_{t+1}-\x_t)= -\J(\x_t)^{\top}\F(\x_t).
\end{align}
Using the Lipschitz continuity of the Gram matrix, we have
\begin{align*}
    &\left\|\G(\x_t)(\x_{t+1}-\x_t)+\J(\x_t)^{\top}\F(\x_t)\right\|\\
    &~~~\leq \left\| \G(\z_t)(\x_{t+1}-\x_t)+\J(\x_t)^{\top}\F(\x_t)\right\|+\left\|(\G(\z_t)-\G(\x_t))(\x_{t+1}-\x_t)\right\|\\
    &\overset{\eqref{eq:prop-G},\eqref{eq:re_org}}{\leq} \lambda_t r_t + 2L_1L_2\|\x_t-\z_t\|r_t,
\end{align*}
and
\begin{align*}
    &\Inner{\J(\x_{t})^{\top}\F(\x_t)+\G(\x_t)(\x_{t+1}-\x_t)}{\x_{t+1}-\x_t}\\
    &=\Inner{\J(\x_{t})^{\top}\F(\x_t)+\G(\z_t)(\x_{t+1}-\x_t)}{\x_{t+1}-\x_t} + \Inner{(\G(\x_t)-\G(\z_t))(\x_{t+1}-\x_t)}{\x_{t+1}-\x_t}\\
    &\!\!\!\!\!\overset{\eqref{eq:prop-G},\eqref{eq:reform}}{\leq}-\lambda_t r_t^2+2L_1L_2r_t^2\|\x_{t}-\z_t\|.
\end{align*}
\end{proof}

\section{The Proof of Lemma~\ref{lm:descent}}\label{apendix:descent}
For the following proof, we denote
\begin{align*}
    \xi_t\triangleq \frac{r_t^2\lambda_t}{6}. 
\end{align*}

\label{sec:proof_descent_lemma}
\begin{proof}
It holds that
\begin{align}
\begin{split}
\label{eq:start}
    &\|\F(\x_{t+1})\|^2\\
    &\overset{\eqref{eq:cubic-growth}}{\leq} \|\F(\x_t)+\J(\x_t)(\x_{t+1}-\x_t)\|^2+M\|\x_{t+1}-\x_t\|^3\\
    &=\|\F(\x_t)\|^2+2\inner{\J(\x_t)^{\top}\F(\x_t)+\J(\x_t)^{\top}\J(\x_t)(\x_{t+1}-\x_t)}{\x_{t+1}-\x_t} \\
    &~~~~~- \|\J(\x_t)(\x_{t+1}-\x_t))\|^2 +M\|\x_{t+1}-\x_t\|^3\\
    &\overset{\eqref{eq:lminner}}{\leq} \|\F(\x_t)\|^2 -2\lambda_tr_t^2 + 4L_2L_1r_t^2\|\x_t-\z_t\|+ Mr_t^3. 
\end{split}
\end{align}
Using Young inequality, we have
\begin{align}
\label{eq:young}
\begin{split}
4L_2L_1\|\x_t-\z_t\|r_t^2 &= \left({4L_2L_1\|\x_t-\z_t\|/c^{2/3}}\right)\cdot c^{2/3}r_t^2 \\
& \leq \frac{1}{3}\cdot\left(4L_2L_1\|\x_t-\z_t\|/c^{2/3}\right)^{3} + \frac{2}{3}\cdot\left(c^{2/3}r_t^2\right)^{3/2}\\
&= \frac{64L_1^{3}L_2^{3}}{3c^2}\|\x_t-\z_t\|^3 + \frac{2c}{3}r_t^3.
\end{split}
\end{align}
Thus, we have
\begin{align*}
    \|\F(\x_{t})\|^2-\|\F(\x_{t+1})\|^2 
    &
\overset{\eqref{eq:start}}{\geq}
      2\lambda_tr_t^2 -Mr_t^3 - 4L_1L_2r_t^2\|\x_t-\z_t\| \\
& = \frac{1}{6}\lambda_t r_t^2 + \frac{11}{6}\lambda_t r_t^2 -  Mr_t^3  - 4L_1L_2r_t^2\|\x_t-\z_t\| \\
& \overset{\eqref{eq:young}}{\geq} 
\frac{1}{6}\lambda_t r_t^2 + \frac{11}{6}\lambda_t r_t^2 -  Mr_t^3 - \left(\frac{64L_1^{3}L_2^{3}}{3c^2}\|\x_t-\z_t\|^3 + \frac{2c}{3}r_t^3\right)\\
&\overset{\eqref{eq:boundrlambda}}{\geq} 
\frac{1}{6}\lambda_t r_t^2 + \frac{11}{6}cr_t^3 -  Mr_t^3 - \left(\frac{64L_1^{3}L_2^{3}}{3c^2}\|\x_t-\z_t\|^3 + \frac{2c}{3}r_t^3\right)\\
&=\frac{1}{6}\lambda_tr_t^2 + \left(\frac{11c}{6}-\frac{2c}{3}-M\right)r_t^3 - \left(\frac{64L_1^{3}L_2^{3}}{3c^2}\|\x_t-\z_t\|^3\right)\\
    &\geq \frac{1}{6}\lambda_tr_t^2 +\frac{c}{6}r_t^3 - \frac{64L_1^3L_2^3}{3c^2}\|\x_t-\z_t\|^3,
\end{align*}
where the last inequality is due to $c\geq M$.
Summing up the above inequality over $t=0,\dots,m-1$, we have
\begin{align*}
    \|\F(\x_{0})\|^2-\|\F(\x_m)\|^2 
    &\geq \sum_{t=0}^{m-1}\xi_t + \frac{c}{6}\sum_{t=0}^{m-1}r_t^3 - \frac{64L_1^3L_2^3}{3c^2}\sum_{t=0}^{m-1}\left(\sum_{i=0}^tr_i\right)^3 \\
    &\geq \sum_{t=0}^{m-1}\xi_t + \left(\frac{c}{6}-\frac{64L_1^3L_2^3m^3}{9c^2}\right)\sum_{t=0}^{m-1}r_t^3\\
    &\geq \sum_{t=0}^{m-1}\xi_t,
\end{align*}
where the first inequality is due to the fact that $\|\x_t-\z_t\|=\|\x_t-\x_0\|\leq \sum_{i=0}^{t-1}\|\x_i-\x_{i+1}\|$, the second inequality is due to  Lemma~\ref{lm:addi-2}, and the last inequality is due to $c\geq 4L_1L_2m$.

\end{proof}

\section{The Proof of Theorem~\ref{thm:LMglobal}}
\begin{proof}
We denote $ \xi_t\triangleq {r_t^2\lambda_t}/{6}. $
Without loss of generality, we suppose $T\equiv 0~({\rm mod}~m)$ and write $T=Km$ for some integer~$K$. 
Lemma~\ref{lm:descent} means we have
\begin{align}
\label{eq:xibound}
    \!\!\Delta_0\defeq  \|\F(\x_0)\|^2\geq \|\F(\x_0)\|^2-\|\F(\x_T)\|^2 \geq \sum_{t=0}^{T-1}\xi_{t}.
\end{align}
Lemma~\ref{lm:pre} means we can lower bound $\xi_t$ as
\begin{align}
\label{eq:xiinequal}
    \xi_t = \frac{r_t^2\lambda_t}{6}\geq \frac{\lambda_t^5}{6c^2(L_1^2+\lambda_t)^2}.
\end{align}
Then we divide the index set $\fI=\{0,1,\cdots, T-1\}$ into 
\begin{align}
\label{eq:setdef}
\begin{split}
   \fI_1\defeq\{i: \lambda_i< L_1^2, i\in\fI\}   ~~~\text{and}~~~\fI_2\defeq\{i: \lambda_i\geq L_1^2,i\in\fI\}.
\end{split}
\end{align}
According to inequality \eqref{eq:xibound}, we have
\begin{align}
\label{eq:Delta_0}
\begin{split}
    \Delta_0 &\overset{\eqref{eq:xibound}}{\geq} \sum_{t=0}^{T-1}\xi_t\\
             &\overset{\eqref{eq:xiinequal}}{\geq} \sum_{t=0}^{T-1}\frac{\lambda_t^5}{6c^2(L_1^2+\lambda_t)^2}\\
             & \!= \sum_{t\in \fI_1} \frac{\lambda_t^5}{6c^2(L_1^2+\lambda_t)^2} + \sum_{t\in\fI_2 }\frac{\lambda_t^5}{6c^2(L_1^2+\lambda_t)^2} \\
             &\overset{\eqref{eq:setdef}}{\geq} \sum_{t\in \fI_1} \frac{\lambda_t^5}{24c^2L_1^2} + \sum_{t\in\fI_2}\frac{\lambda_t^3}{24c^2}\\
             &\overset{\eqref{eq:setdef}}{\geq}  \sum_{t\in \fI_1} \frac{\lambda_t^5}{24c^2L_1^2} + \sum_{t\in\fI_2}\frac{L_1^6}{24c^2}.
\end{split}
\end{align}
This indicates that the cardinality $\fI_2$ satisfies 
\begin{align*}
     \frac{|\fI_2|L_1^6}{24c^2}\leq \Delta_0.
\end{align*}
Thus, taking the total iteration number
\begin{align*}
       T=\left\lceil\frac{24\Delta_0 c^2}{L_1^6} +\frac{24\Delta_0L_1^2}{c^{0.5}\epsilon^{2.5}}\right\rceil
\end{align*}
leads to the cardinality of $\fI_1$ satisfies 
\begin{align}
\label{eq:fI1}
\begin{split}
    |\fI_1|&=T-|\fI_2|\geq \frac{24\Delta_0 c^2}{L_1^6} +\frac{24\Delta_0L_1^2}{c^{0.5}\epsilon^{2.5}} - \frac{24\Delta_0 c^2}{L_1^6} = \frac{24\Delta_0L_1^2}{c^{0.5}\epsilon^{2.5}}.
\end{split}
\end{align}
Then we have
\begin{align*}
\min_{t\in\{0,\cdots T-1\}} \|\J(\x_t)^{\top}\F(\x_t)\|\leq \min_{t\in\fI_1}\|\J(\x_t)^{\top}\F(\x_t)\| \overset{\eqref{eq:Delta_0}}{\leq} \left(\frac{24\Delta_0 L_1^2}{c^{0.5}|\fI_1|}\right)^{2/5} \overset{\eqref{eq:fI1}}{=} \epsilon,
\end{align*}
where we take $c = \max\{4L_1L_2m, M\}$ by following the setting of Lemma~\ref{lm:descent}.
Therefore, we finish the proof.
\end{proof}

\section{The Proof of Lemma~\ref{lm:superlm}}
\label{sec:proof_local_lemma}
\begin{proof}
Proposition \ref{prop:local} implies
    \begin{align}        \label{eq:Gregubound} &\Norm{(\G(\z_t)+\lambda_t\I)^{-1}}\leq\Norm{(\G(\z_t))^{-1}}\leq \frac{2}{\mu^2},
\end{align}
\begin{align*}
   \lambda_t &= \sqrt{c\|\J(\x_t)^{\top}\F(\x_t)\|}\\
   &{\leq} \sqrt{cL_1\|\F(\x_t)\|} \nonumber\\
        &\overset{\eqref{eq:unique}}{=}c^{0.5}L_1^{0.5}\|\F(\x_t)-\F(\x^*)\|^{0.5}\\
        &\overset{\eqref{eq:lipF}}{\leq} c^{0.5}L_1\|\x_t-\x^*\|^{0.5}, 
    \end{align*}
    and
    \begin{align}
    \label{eq:boundFJ}
    \begin{split}
       &\| \J(\x_t)(\x_t-\x^*)-\F(\x_t)\|
       \\&\overset{\eqref{eq:unique}}{=} \|\J(\x_t)(\x_t-\x^*) - (\F(\x_t)-\F(\x^*))\|\\
       &
       \overset{\eqref{eq:lipF}}{\leq} \frac{L_2}{2}\|\x_t-\x^*\|^2.
    \end{split}
    \end{align}
    The update rule of Algorithm~\ref{alg:LLM} means 
    \begin{align*}
   \begin{split}
    &\x_{t+1}-\x^* \\
    &= \x_t-\x^* -\left(\G(\z_t)+\lambda_t\I\right)^{-1}\J(\x_t)^{\top}\F(\x_t)\\
    &= \left(\G(\z_t)+\lambda_t\I\right)^{-1}\left((\G(\z_t)+\lambda_t\I)(\x_t-\x^*)-\J(\x_t)^{\top}\F(\x_t)\right)\\
    & = \left(\G(\z_t)+\lambda_t\I\right)^{-1} \left(\J(\x_t)^{\top}(\J(\x_t)(\x_t-\x^*)-\F(\x_t))+ (\G(\z_t)-\G(\x_t))(\x_t-\x^*)+\lambda_t(\x_t-\x^*))\right),
    \end{split}
    \end{align*}
    then it holds that
    \begin{align*}
        &\|\x_{t+1}-\x^*\|\\
        & \overset{\eqref{eq:prop-G},\eqref{eq:Gregubound}}{\leq} \frac{2}{\mu^2}\left(\|\J(\x_t)\|\| \J(\x_t)(\x_t-\x^*)-\F(\x_t)\|+ 2L_1L_2\|\z_t-\x_t\|\|\x_t-\x^*\|+\lambda_t\|\x_t-\x^*\|\right)\\
        & ~~\overset{\eqref{eq:boundFJ}}{\leq}
        \frac{L_1L_2}{\mu^2} \|\x_t-\x^*\|^2 + \frac{2L_1L_2}{\mu^2}\|\x_t-\x^*\|\|\z_t-\z^*\| + \frac{2}{\mu^2}\lambda_t\|\x_t-\x^*\|\\
        &~~\overset{\eqref{eq:boundrlambda}}{\leq} \frac{L_1L_2}{\mu^2}\|\x_t-\x^*\|^2 + \frac{2L_1c^{0.5}}{\mu^{2}}\|\x_t-\x^*\|^{1.5} + \frac{2L_1L_2}{\mu^2}\|\z_t-\x^*\|\|\x_t-\x^*\|\\
        &~~~=\alpha_1 \|\x_t-\x^*\|^{2} +\alpha_2 \|\x_t-\x^*\|^{1.5}+\beta \|\z_t-\x^*\|\|\x_t-\x^*\|.
    \end{align*}
    Since we suppose $\x_t,\z_t\in\fS$, it holds that
    \begin{align*}
    &\|\x_{t+1}-\x^*\|\\
    &=
         (\alpha_1 \|\x_t-\x^*\| +\alpha_2 \|\x_t-\x^*\|^{0.5}+\beta \|\z_t-\x^*\|)\|\x_t-\x^*\|\\
         &\overset{\eqref{eq:local-condi}}{\leq} \|\x_t-\x^*\|,
    \end{align*}
    which means $\x_{t+1}\in\fS$.
\end{proof}
\section{The Proof of Theorem~\ref{thm:LMlocal}}
\begin{proof}
We first use induction to prove that all points $\x_0,\dots,\x_T$ generated by Algorithm~\ref{alg:LLM} are in the set~$\fS$ defined in Lemma \ref{lm:superlm}:
\begin{itemize}[leftmargin=0.5cm]
    \item The initial condition \eqref{eq:initial} means $\x_0\in\fS$.
    \item We suppose $\vx_1,\dots,\vx_{t}\in\fS$ for any $t=1,\dots,T-1$.
    Then we have $\z_{t}\in\{\vx_1,\dots,\vx_{t}\}\subseteq\fS$, which implies $\vz_t\in\fS$. According to Lemma~\ref{lm:superlm}, we have $\x_{t+1}\in\fS$. This finishes the induction.    
\end{itemize}  
Since points $\x_0,\cdots,\x_T$ are all in the set $\fS$, Lemma \ref{lm:superlm} means \eqref{eq:iterlocal} holds for all $t=0,\dots,T$. We let $\eta = 2(\alpha_1+\alpha_2^2)$ and denote $s_t \triangleq \eta\|\x_t-\x^*\|$, then it holds that
\begin{align*}
    s_{t+1}\overset{\eqref{eq:iterlocal}}{\leq} \frac{1}{2} \left(s_t^2+s_t^{1.5}\right) + s_{\pi(t)}s_t.
\end{align*}
Since it holds $s_0\overset{\eqref{eq:initial}}{\leq} {1}/{2^{4}}$, we achieve
\begin{align*}
    s_{t}\leq \left(\frac{1}{2}\right)^{2(1+(1+m/2)^{\pi(t)})(1+(t\%m)/2)}
\end{align*}
by using the results of Theorem D.3 of \citet{doikov2023second} (see Lemma~\ref{lm:addi-2} in Appendix~\ref{sec:usefulsequece}), which finishes the proof.
\end{proof}

\section{Baseline Algorithms in Experiment}
\label{sec:appbase}
We present the two baselines GD and LM in Algorithm~\ref{alg:GD} and \ref{alg:LM} respectively.

\begin{algorithm}[t]
\caption{Gradient Descent for Solving Nonlinear Equations (GD)}\label{alg:GD}
\begin{algorithmic}[1]
\STATE \textbf{Input:} $\x_0$, $\eta$, and $T$ \\[0.15cm]
\STATE \textbf{for} $t=0,1,\dots T-1$ \\[0.15cm]
\STATE \quad    $\x_{t+1}=\x_t-\eta\J(\x_t)^{\top}\F(\x_t)$\\[0.15cm]
\STATE \textbf{end for}\\[0.15cm]
\end{algorithmic}
\end{algorithm}

\begin{algorithm}[t]
\caption{Levernberg--Marquardt~\cite{mishchenko2021regularized}  (LM)}\label{alg:LM}
\begin{algorithmic}[1]
\STATE \textbf{Input:} $\x_0$, $c$, and $T$ \\[0.15cm]
\STATE \textbf{for} $t=0,1,\dots T-1$ \\[0.15cm]
\STATE \quad    $\x_{t+1}=\x_t-\left(\J(\x_t)^{\top}\J(\x_t)+\sqrt{c\|\J(\x_t)^{\top}\F(\x_t)\|}\I\right)^{-1}\J(\x_t)^{\top}\F(\x_t)$\\[0.15cm]
\STATE \textbf{end for}\\[0.15cm]
\end{algorithmic}
\end{algorithm}

\end{document}